\documentclass[11pt,leqno,a4paper]{article}

\usepackage[hypertex]{hyperref}
\usepackage{amsfonts}
\usepackage{latexsym}
\usepackage{amssymb,amsmath}
\usepackage{epsfig}
\newtheorem{theorem}{Theorem}[section]
\newtheorem{corollary}[theorem]{Corollary}
\newtheorem{lemma}[theorem]{Lemma}

\newtheorem{definition}[theorem]{Definition}
\newtheorem{remark}[theorem]{Remark}
\newcommand{\ep}{\varepsilon}
\newcommand{\eps}[1]{{#1}_{\varepsilon}}
\newcommand{\R}{\mathbb R}

\newcommand{\N}{\mathbb N}

\def\qd{{\rule{2.3mm}{2.3mm}}}
\def\qed{{\hfill$\quad$\qd\\}}

\date{February 3, 2010}
\begin{document}

%% Place the running title of the paper with 40 letters or less in []
 %% and the full title of the paper in { }.
\title{Uniform stabilization in weighted Sobolev spaces 
for the KdV equation posed on the half-line}

% Place all authors' names in [ ] shown as running head;
% No more than 40 letters. Leave { } empty
% Please use `and' to connect the last two names if applicable
\author{Ademir F. Pazoto 
\thanks{Instituto de Matem\'atica,
Universidade Federal do Rio de Janeiro,
P.O. Box 68530, CEP 21945-970, Rio de Janeiro, RJ,  Brasil
({\tt ademir@im.ufrj.br})}
\and Lionel Rosier
\thanks{Institut Elie Cartan, UMR 7502 UHP/CNRS/INRIA,
B.P. 239, F-54506 Vand\oe uvre-l\`es-Nancy Cedex, France
({\tt rosier@iecn.u-nancy.fr})} }

\maketitle

%The abstract of your paper
\begin{abstract}
Studied here is the large-time behavior of solutions of the Korteweg-de Vries
equation posed on the right half-line under the effect of a localized damping.
Assuming as in \cite{linares-pazoto} that the damping is active on a set
$(a_0,+\infty)$ with $a_0>0$, we establish the exponential 
decay of the solutions in the weighted spaces $L^2((x+1)^mdx)$ for $m\in \N ^*$
and $L^2(e^{2bx}dx)$ for $b>0$ by a Lyapunov approach. The decay of the spatial derivatives
of the solution is also derived.

\end{abstract}

%The title of your section 1
\section{Introduction}
The Korteweg-de Vries (KdV) equation was first derived as a model
for the propagation of small amplitude long water waves along a channel
\cite{bouss,jager,korteweg}. It has been intensively studied from
various aspects for both mathematics and physics since the 1960s when
solitons were discovered through solving the KdV equation, and the
inverse scattering method, a so-called nonlinear Fourier transform,
was invented to seek solitons \cite{gardner,miura}. It is now
well known that the KdV equation is not only a good model for water
waves but also a very useful approximation model in nonlinear
studies whenever one wishes to include and balance weak nonlinear
and dispersive effects. 

The initial boundary value problems (IBVP) arise
naturally in modeling small-amplitude long waves in a channel with a
wavemaker mounted at one end \cite{bona1,bona2,bona3,R04}. Such
mathematical formulations have received considerable attention in
the past, and a satisfactory theory of global well-posedness is available for initial and boundary conditions satisfying physically
relevant smoothness and consistency assumptions (see e.g.
\cite{bona1,bona4,bona6,bona7,colliander,faminskii,faminskii2}  
and the references therein).

The analysis of the long-time behavior of IBVP on the 
quarter-plane for KdV has also received considerable attention over recent years,
and a review of some of the results related to the issues we address
here can be found in \cite{bona5,bona7,leach}. For stabilization
and controllability issues on the half line, we refer the reader to  \cite{linares-pazoto} and \cite{R00,R02}, respectively. 

In this work, we are concerned with the asymptotic behavior of the
solutions of the IBVP for the KdV
equation posed on the positive half line under the presence of a
localized damping represented by the function $a$; that is,
\begin{equation} \label{1}
\begin{cases}
u_t + u_x + u_{xxx} + uu_x + a(x)u = 0, \quad x,\, t \in \mathbb{R}^+,\\
u(0,t)=0, \quad t>0,\\
u(x,0)=u_0(x), \quad x > 0.
\end{cases}
\end{equation}

%In \cite{bona1} Bona and Winther proposed the above boundary-value
%problem ($a\equiv 0$) to describe the evolution of unidirectional
%waves generated at one end of a homogeneous stretch of a certain
%medium and which are allowed to  propagate into the initially
%undisturbed medium beyond a wavemaker. They also gave the first
%result regarding well-posedness for the IBVP \eqref{1}.

Assuming $a(x)\ge 0$ a.e. and that $u(.,t)\in H^3(\R ^+ )$, it follows from a simple computation that
\begin{eqnarray}\label{3}
\frac{dE}{dt} = - \int_0^\infty a(x)|u(x,t)|^2 dx  - \frac{1}{2}
|u_x(0,t)|^2
\end{eqnarray}
where
\begin{eqnarray}\label{4}
E(t) = \frac{1}{2}\,\int_0^\infty |u(x,t)|^2 dx
\end{eqnarray}
is the total energy associated with (\ref{1}). Then, we see that the term $a(x)u$ plays the role of a feedback damping mechanism and, consequently, it is
natural to wonder whether the solutions of (\ref{1}) tend to zero as
$t\rightarrow \infty$ and under what rate they decay. When
$a(x)>a_0>0$ almost everywhere in $\mathbb{R}^+$, it is very simple
to prove that $E(t)$ converges to zero as $t$ tends to infinity. The
problem of stabilization when the damping is effective only in a
subset of the domain is much more subtle. The following result was 
obtained in \cite{linares-pazoto}. 
\begin{theorem}
\label{thm0} 
Assume that the function $a=a(x)$ satisfies the following property
\begin{equation}
\label{2}
a\in L^\infty(\R ^+),\ a\ge 0\mbox{ a.e. in }\ \R ^+ \mbox{ and } 
a(x)\ge a_0>0 \mbox{ a.e. in } \ (x_0,+\infty) 
\end{equation}
for some numbers $a_0,x_0>0$. Then for all $R>0$ there exist two numbers
$C>0$ and $\nu >0$ such that for all $u_0\in L^2(\R ^+)$ with 
$||u_0||_{L^2(\R ^+ )}\le R$, the solution $u$ of \eqref{1} satisfies
\begin{equation}
\label{L2}
||u(t)||_{L^2(\R ^+)} \le C e^{-\nu t} 
||u_0||_{L^2(\R ^+)}\cdot
\end{equation}
\end{theorem}
Actually, Theorem \ref{thm0} was proved in \cite{linares-pazoto} under the 
additional hypothesis that 
\begin{equation}
\label{L3}
a(x)\ge a_0 \ \mbox{ a.e. in }\  (0,\delta ) 
\end{equation}
for some $\delta >0$, but \eqref{L3} may be dropped by replacing the unique continuation property 
\cite[Lemma 2.4]{linares-pazoto} by \cite[Theorem 1.6]{rosier-zhang}.
The exponential decay of $E(t)$ is obtained following the methods in \cite{pazoto,PMVZ,R97} which
combine multiplier techniques and compactness arguments to reduce
the problem to some unique continuation property for weak
solutions of KdV. 

Along this work we assume that the real-valued function $a=a(x)$
satisfies the condition \eqref{2} for some given positive numbers $a_0,x_0$. 
%\begin{equation}\label{2}
%\vspace{3mm}a\in W^{3,\infty}(\R ^+),\ \ a(x)\geq 0\ \hbox{ a.e. in }\ \R^+,\
%a(x)\geq a_0 > 0 \  \hbox { a.e. in } \ \omega :=(x_0,+\infty),
%\end{equation}
%where $x_0 > 0$ is a given number. 
In this paper we investigate the stability properties 
of \eqref{1} in the weighted spaces introduced by Kato in \cite{kato}. 
More precisely, for $b> 0$ and $m\in \mathbb{N}$, we prove
that the solution $u$ exponentially decays to $0$ in $L^2_b$ and
$L^2_{(x+1)^m dx}$ (if $u(0)$ belongs to one of these spaces), where
$$L^2_b = \{u:\mathbb{R}^+\rightarrow \mathbb{R} ; \int_0^\infty |u(x)|^2 e^{2bx} dx <
\infty \},$$
$$L^2_{(x+1)^m dx} = \{u:\mathbb{R}^+\rightarrow \mathbb{R}; \int_0^\infty |u(x)|^2 (x + 1)^m dx < \infty \} .$$
The following weighted Sobolev spaces 
$$H^s_b = \{u:\mathbb{R}^+\rightarrow \mathbb{R};\  \partial _x
^i u\in L^2_b \ \mbox{ for } 0 \le i\le s;\ 
u(0)=0 \hbox{  if } s\ge 1 \}$$
and
$$H^s_{(x+1)^mdx} = \{u:\mathbb{R}^+\rightarrow \mathbb{R} ; 
\,\partial^i_x u \in L^2_{(x+1)^{m-i} dx}\,\hbox{ for }\, 0\leq i \leq s ;\ 
u(0)=0 \hbox{  if } s\ge 1\},$$
endowed with their usual inner products, will be used thereafter. 
Note that $H^0_b=L^2_b$ and that $H^0_{(x+1)^mdx}=L^2_{(x+1)^mdx}$.

The exponential decay in $L^2_{(x+1)^m dx}$ is obtained by
constructing a convenient Lyapunov function (which actually decreases
strictly on the sequence of times $\{ kT \}_{k\ge 0}$) by induction on $m$.  
For $u_0\in L^2_{(x+1)^m dx}$, we also prove the following estimate
\begin{equation}
||u(t)||_{H^1_{(x+1)^mdx}} 
\le C\frac{e^{-\mu t}}{\sqrt{t}} ||u_0||_{L^2_{(x+1)^mdx}}
\label{globalkato}
\end{equation}
in two situations:
(i)  $m=1$    and  $||u_0||_{L^2_{(x+1)^m dx}}$ is arbitrarily large;
(ii) $m\ge 2$ and  $||u_0||_{L^2_{(x+1)^m dx}}$ is small enough.
In the situation (ii), we first establish a similar estimate for the linearized 
system and next apply the contraction mapping principle in a space of functions 
fulfilling the exponential decay.
Note that \eqref{globalkato} combines the (global) Kato smoothing effect
to the exponential decay. 

The exponential decay in $L^2_b$ is established for any initial data 
$u_0\in L^2_b$ under the additional assumption that $4 b^3 + b < a_0$. 
Next, we can derive estimates of the form
$$
||u(t)||_{H^s_b} \le C \frac{e^{-\mu t}}{t^{s/2}} ||u_0||_{L^2_b}
$$ 
for any $s\ge 1$, revealing that $u(t)$ decays exponentially to 0 in 
strong norms. 

It would be interesting to see if such results are still true when the
function $a$ has a smaller support. It seems reasonable to conjecture that
similar positive results can be derived when the support of $a$ contains 
a set of the form $\cup_{k\ge 1}[ka_0,ka_0+b_0]$ where $0<b_0<a_0$, while
a negative result probably holds when the support of $a$ is a finite interval,
as the $L^2$ norm of a soliton-like initial data may not be sufficiently dissipated over time. Such issues will be discussed elsewhere.

The plan of this paper is as follows. Section 2 is devoted to 
global well-posedness results in the weighted spaces $L^2_b$ and 
$L^2_{(x+1)^2dx}$. In section 3,
we prove the exponential decay in $L^2_{(x+1)^mdx}$ and $L^2_b$, and establish
the exponential decay of the derivatives as well.

\section{Global well-posedness}

\subsection{Global well-posedness in $\bf L^2_b$}
Fix any $b>0$.
To begin with, we apply the classical semigroup theory to the
linearized system
\begin{equation} \label{linear}
\begin{cases}
u_t + u_x + u_{xxx} + a(x)u = 0, \quad x,\, t \in \mathbb{R}^+,\\
u(0,t)=0, \quad t>0,\\
u(x,0)=u_0(x), \quad x > 0.
\end{cases}
\end{equation}
Let us consider the operator
$$A: D(A)\subset L^2_b \rightarrow L^2_b$$
with domain
$$D(A) = \{ u\in L^2_b;\ \partial_x^i u\in L^2_b \mbox{ for } 1\leq i \leq 3\,\,\mbox{and}\,\,u(0)=0\}$$
defined by
$$Au = - u_{xxx} - u_x - a(x)u.$$
Then, the following result holds.

\begin{lemma}\label{semigroup} The operator $A$ defined above generates a continuous
semigroup of operators $(S(t))_{t\geq 0}$ in $L^2_b$.
\end{lemma}

\noindent {\bf Proof.} We first introduce the new variable $v=e^{bx}u$ and
consider the following (IBVP)
\begin{equation} \label{1.1}
\begin{cases}
v_t + (\partial_x -b)v + (\partial_x - b)^3v  + a(x)v = 0, \quad x,\, t \in \mathbb{R}^+,\\
v(0,t)=0, \quad t>0,\\
v(x,0)=v_0(x)=e^{bx}u_0(x), \quad x > 0.
\end{cases}
\end{equation}
Clearly, the operator $B: D(B)\subset L^2(\mathbb{R}
^+) \rightarrow
L^2(\mathbb{R}^+)$ with domain
$$D(B) = \{u\in H^3(\mathbb{R}^+);\ u(0)=0\}$$
defined by
$$Bv = - (\partial_x -b)v - (\partial_x - b)^3v  - a(x)v$$
is densely defined and closed. So, we are done if we prove that
for some real number $\lambda$ the operator $B-\lambda $
and its adjoint $B^\ast -\lambda $ are both dissipative in $L^2(\mathbb{R}^+)$.
It is readily seen that $B^\ast : D(B^\ast)\subset
L^2(\mathbb{R}^+)\rightarrow L^2(\mathbb{R}^+)$ is given by $B^\ast v= (\partial_x
+b)v + (\partial_x + b)^3v - a(x)v$ with domain
$$D(B^\ast) = \{v\in H^3(\mathbb{R}^+);\ v(0)=v'(0)=0\}.$$
Pick any $v\in D(B)$. After some integration by parts, we obtain
that
$$(Bv,v)_{L^2} = -\frac{1}{2}v_x^2(0) - 3b\int_0^\infty v^2_x dx + (b+b^3)\int_0^\infty v^2dx
- \int_0^\infty a(x)v^2 dx,$$ that is,
$$([B - (b^3 + b) ]v,v)_{L^2} \leq 0.$$
Analogously, we deduce that for any $v\in D(B^\ast)$
$$(v,[B^\ast - (b^3 + b) ]v)_{L^2} \leq 0$$
which completes the proof.
\qed

The following linear estimates will be needed.

\begin{lemma} \label{lin} Let $u_0 \in L^2_b$ and $u = S(\cdot)u_0$. Then, for any
$T>0$
\begin{equation}\label{en-id1}
\frac{1}{2}\int_0^\infty |u(x,T)|^2 dx - \frac{1}{2}\int_0^\infty
|u_0(x)|^2 dx +\int_0^T\int_0^\infty a(x)|u|^2 dxdt
+\frac{1}{2}\int_0^T u_x^2(0,t)dt = 0
\end{equation}
\begin{equation}\label{en-id2}
\begin{array}{l}
\displaystyle\frac{1}{2}\int_0^\infty |u(x,T)|^2 e^{2bx}dx -
\frac{1}{2}\int_0^\infty |u_0(x)|^2 e^{2bx}dx
+\,3b\displaystyle\int_0^T\int_0^\infty u^2_x e^{2bx} dx dt\\
- \displaystyle(4b^3 + b)\int_0^T\int_0^\infty u^2 e^{2bx} dx dt
+\displaystyle\int_0^T\int_0^\infty a(x)|u|^2 e^{2bx}dxdt
+\frac{1}{2}\int_0^T u_x^2(0,t)dt = 0.
\end{array}
\end{equation}
As a consequence,
\begin{equation}\label{lin-est}
\begin{array}{l}
||u||_{L^\infty (0,T;L^2_b)} + ||u_x||_{L^2(0,T;L^2_b)}\leq
C\,||u_0||_{L^2_b},
\end{array}
\end{equation}
where $C=C(T)$ is a positive constant.
\end{lemma}

\noindent {\bf Proof.}
Pick any $u_0\in D(A)$. Multiplying the equation in
(\ref{1}) by $u$ and integrating over $(0,+\infty )\times(0,T)$, we obtain
(\ref{en-id1}). Then, the identity may be extended to any initial
state $u_0\in L^2_b$ by a density argument. To derive (\ref{en-id2})
we first multiply the equation by $(e^{2bx} - 1)u$ and integrate by
parts over $(0,+\infty )\times (0,T)$ to deduce that
\begin{equation}
\begin{array}{l}
\vspace{1mm}\displaystyle\frac{1}{2}\int_0^\infty |u(x,T)|^2(e^{2bx}-1)dx -
\frac{1}{2}\int_0^\infty |u_0(x)|^2 (e^{2bx}-1)dx \,+\\
\vspace{1mm}+3b\displaystyle\int_0^T\int_0^\infty u^2_x e^{2bx} dx dt
- \displaystyle(4b^3 + b)\int_0^T\int_0^\infty u^2 e^{2bx} dx dt \, + \\
+\displaystyle\int_0^T\int_0^\infty a(x)|u|^2 (e^{2bx}- 1) dxdt = 0.
\nonumber
\end{array}
\end{equation}
Adding the above equality and (\ref{en-id1}) hand to hand, we
obtain (\ref{en-id2}) using the same density argument.
Then, Gronwall inequality, (\ref{2}) and (\ref{en-id2}) imply that
$$||u||_{L^\infty(0,T;L^2_b)} \leq \, C\,||u_0||_{L^2_b},$$
with $C=C(T) > 0$. This estimate together with (\ref{en-id2})
%and Gronwall inequality
gives us
$$||u_x||_{L^2(0,T;L^2_b)}\leq C\,||u_0||_{L^2_b},$$
where $C=C(T)$ is a positive constant.
\qed

The global well-posedness result reads as follows:

\begin{theorem}
\label{global-exp} For any $u_0 \in L^2_b$ and any $T>0$, there exists a unique solution
$u\in C([0,T];L^2_b)\cap L^2(0,T;H^1_b)$ of \eqref{1}.
\end{theorem}

\noindent {\bf Proof.}
By computations similar to those performed in the proof of Lemma
\ref{lin}, we obtain that for any $f\in C^1([0,T];L^2_b)$ and any
$u_0\in D(A)$, the
solution $u$ of the system
$$
\left\{
\begin{array}{ll}
u _t + u_x + u_{xxx} + a(x)u =f, \quad & x\in \mathbb{R}^+,\ t\in (0,T), \\
u(0,t)=0, & t\in (0,T),\\
u(x,0)=u_0(x),& x\in \mathbb{R}^+,
\end{array}
\right.
$$
fulfills
\begin{equation}\label{est-ex}
\sup_{0\le t\le T} ||u(t)||_{L^2_b} +(\int_0^T\!\!\!\int_0^\infty
|u_x|^2 e^{2bx}dxdt)^{\frac{1}{2}} \leq C\left( ||u_0||_{L^2_b} +\int_0^T
||f||_{L^2_b} dt \right)
\end{equation}
for some constant $C=C(T)$ nondecreasing in $T$. A density argument
yields that $u\in C([0,T]; L^2_b)$ when $f\in L^1(0,T;L^2_b)$
and $u_0\in L^2_b$.

Let $u_0\in L^2_b$ be given. To prove the existence of a solution of
\eqref{1} we introduce the map $\Gamma$ defined by
$$
(\Gamma u)(t)=S(t)u_0+\int_0^tS(t-s)N(u(s))\, ds
$$
where $N(u)= -uu_x$, and the space
$$F = C([0,T];L^2_b)\cap L^2(0,T;H^1_b)$$
endowed with its natural norm.
We shall prove that $\Gamma$ has a fixed-point in some ball $B_R(0)$
of $F$. We need the following \\

{\sc Claim 1.} If $u\in H^1_b$ then
$$||u^2e^{2bx}||_{L^\infty (\R ^+)} \leq (2 + 2b)\,
||u||_{L^2_b} ||u||_{H^1_b}.$$

From Cauchy-Schwarz inequality, we get for any $\overline{x}\in
\mathbb{R}^+$
$$
\begin{array}{l}
u^2(\overline{x})e^{2b\overline{x}}
= \displaystyle\int_0^{\overline{x}} [u^2e^{2bx}]_x dx
= \int_0^{\overline{x}}
[2uu_xe^{2bx} + 2bu^2e^{2bx}]dx \\
\leq 2(\displaystyle\int_0^\infty u^2
e^{2bx}dx)^{\frac{1}{2}}(\int_0^\infty u_x^2
e^{2bx}dx)^{\frac{1}{2}} + 2b\int_0^\infty u^2 e^{2bx} dx \leq (2 +
2b) ||u||_{L^2_b} ||u||_{H^1_b}
\end{array}
$$
which guarantees that Claim 1 holds.

{\sc Claim 2.} There exists a constant $K>0$ such that for $0<T\le 1$
\begin{equation}
||\Gamma(u) - \Gamma(v) ||_F \leq KT^{\frac{1}{4}} (||u||_F +
||v||_F) ||u - v||_F, \quad \forall\,u, v\, \in F.\nonumber
\end{equation}

According to the previous analysis,
\begin{equation}\label{fp-0}
||\Gamma(u) - \Gamma(v) ||_F \leq C ||uu_x -
vv_x||_{L^1(0,T;L^2_b)}.\nonumber
\end{equation}
So, applying triangular inequality and H\"older inequality, we have
\begin{eqnarray}
&&||\Gamma(u) - \Gamma(v) ||_F \leq C \{||u -
v||_{L^2(0,T;L^\infty(0,\infty))}||u||_{L^2(0,T;H^1_b)} +
\nonumber\\
&&\qquad\qquad + ||v||_{L^2(0,T;L^\infty(0,\infty))}||u - v||_{L^2(0,T;H^1_b)} \}. \label{fp-1}
\end{eqnarray}
Now, by Claim 1, we have
\begin{equation}\label{fp-2}
\begin{array}{l}
||u||_{L^2(0,T;L^\infty(0,\infty))}\leq C\, T^{\frac{1}{4}}
||u||_{L^\infty(0,T;L^2_b)}^{\frac{1}{2}}||u||_{L^2(0,T;H^1_b)}^{\frac{1}{2}}.
\end{array}
\end{equation}
Then, combining (\ref{fp-1}) and (\ref{fp-2}), we deduce that
\begin{equation}
\label{fp-3}
\begin{array}{l}
||\Gamma(u) - \Gamma(v) ||_F \leq C\, T^{\frac{1}{4}} \{\,||u||_F +
||v||_F\,\}||u - v||_F.
\end{array}
\end{equation}

Let $T>0$, $R>0$ be numbers whose values will be specified later,
and let $u\in B_R(0)\subset F$ be given. Then, by Claim 2 and Lemma
\ref{lin}, $\Gamma u \in F$ and
$$
||\Gamma u||_F \leq C\,(\,||u_0||_{L^2_b} +
T^{\frac{1}{4}}||u||^2_F\,).
$$
Consequently, for $R=2C||u_0||_{L^2_b}$ and $T>0$ small enough, $\Gamma$ maps
$B_R(0)$ into itself. Moreover, we infer from \eqref{fp-3} that this mapping
contracts if $T$ is small enough. Then, by the contraction mapping
theorem, there exists a unique solution $u\in B_R(0)\subset F$ to the
problem \eqref{1} for $T$ small enough.

In order to prove that this solution is global, we need some a
priori estimates. So, we proceed as in the proof of Lemma \ref{lin}
to obtain for the solution $u$ of \eqref{1}
\begin{equation}\label{fp30}
\frac{1}{2}\int_0^\infty |u(x,T)|^2 dx - \frac{1}{2}\int_0^\infty
|u_0(x)|^2 dx +\int_0^T\int_0^\infty a(x)|u|^2 dxdt
+\frac{1}{2}\int_0^T u_x^2(0,t)dt = 0
\end{equation}
and
\begin{eqnarray}
\vspace{1mm}\displaystyle\frac{1}{2}\int_0^\infty |u(x,T)|^2 e^{2bx}dx -
\frac{1}{2}\int_0^\infty |u_0(x)|^2 e^{2bx}dx + \frac{1}{2}\int_0^T u_x^2(0,t)dt
\nonumber\\
\vspace{1mm}+\,3b\displaystyle\int_0^T\int_0^\infty u^2_x e^{2bx} dx dt - \displaystyle(4b^3 + b)\int_0^T\int_0^\infty u^2 e^{2bx} dx dt\nonumber\\
+\displaystyle\int_0^T\int_0^\infty a(x)|u|^2 e^{2bx}dxdt
-\frac{2b}{3}\int_0^T \int_0^\infty u^3 e^{2bx}dx dt= 0.\label{fp-4}
\end{eqnarray}
%From (\ref{fp30}), we have
%$$||u||_{L^\infty(0,T;L^2(0,\infty))} \leq ||u_0||_{L^2_b}.$$
First, observe that
$$|\int_0^\infty u^2 e^{2bx} dx| = |-\frac{1}{b}\int_0^\infty uu_x e^{2bx}dx|
\leq \frac{1}{b} (\int_0^\infty u^2 e^{2bx}
dx)^{\frac{1}{2}}(\int_0^\infty u^2_x e^{2bx} dx)^{\frac{1}{2}},$$
therefore,
\begin{equation}
\begin{array}{l}%\label{exp5.1}
\displaystyle\int_0^\infty u^2 e^{2bx} dx \leq
\frac{1}{b^2}\int_0^\infty u^2_x e^{2bx} dx.\nonumber
\end{array}
\end{equation}
Combined to Claim 1, this yields
$$
||u(x)e^{bx}||_{L^\infty (\R ^+)} 
\le C ||u_x ||_{L^2_b}.
$$
On the other hand, it follows from \eqref{fp30}
that
$$
|| u(t) ||_{L^2(\R ^+)} \le ||u_0||_{L^2(\R ^+)},
$$
hence 
%$$
%\int_0^\infty |u|^3 e^{2bx} dx
%\le ||u||_{L^2} ||ue^{bx}||_{L^2} ||ue^{bx} ||_{L^\infty}
%\le c||u||_{L^2_b}(||u_x||_{L^2_b} + ||u||_{L^2_b})$$
%hence
\begin{eqnarray*}
\int_0^T \!\!\! \int_0^\infty |u|^3 e^{2bx} dx dt 
&\leq& \int_0^T ||ue^{bx}||_{L^\infty(\R ^+)}
(\int_0^\infty |u|^2e^{bx}dx)dt \\
&\le& C \int_0^T ||u_x||_{L^2_b} ||u||_{L^2_b}||u||_{L^2}dt \\
&\le& \delta ||u_x||^2_{L^2(0,T;L^2_b)} +
 C_\delta ||u||^2_{L^2(0,T;L^2_b)}, 
\end{eqnarray*}
where $\delta > 0 $ is arbitrarily chosen
and $C=C(b,\delta, ||u_0||_{L^2(\R ^+)} )$ is a positive constant. 
Combining this inequality (with $\delta <9/2$) to  (\ref{fp-4}) results in 
$$
||u(T)||^2_{L^2_b} 
\le ||u_0||^2_{L^2_b} + C\int_0^T||u||^2_{L^2_b}dt
$$ 
where $C=C(b,||u_0||_{L^2(\R ^+)})$ does not depend on $T$. It follows from 
Gronwall lemma that 
$$
||u(T)||^2_{L^2_b} \le ||u_0||^2_{L^2_b} e^{CT}
$$
for all $T>0$, which gives the global well-posedness.
\qed

\subsection{Global well-posedness in $L^2_{(x+1)^2dx} $}

%Our first result is a consequence of Theorem \ref{global-exp}:

\begin{definition}
For $u_0\in L^2_{(x+1)^2dx}$ and $T>0$, we denote by a {\em mild solution} 
of \eqref{1} any function $u\in C([0,T];L^2_{(x+1)^2 dx})\cap 
L^2(0,T;H^1_{(x+1)^2dx})$ which solves \eqref{1}, and such that for some 
$b>0$ and some sequence $\{ u_{n,0}\}\subset L^2_b$ we have 
\begin{eqnarray*}
&&u_{n,0}\to u_0 \,\, \mbox{ strongly in }\,\,  L^2_{(x+1)^2dx},\\
&&u_n \to u\,\,\mbox{weakly}\ast\,\,\mbox{in}\,\, 
L^\infty (0,T; L^2_{(x+1)^2dx}),\\
&&u_n \to u\,\,\mbox{weakly}\,\,\mbox{in}\,\,
L^2(0,T; H^1_{(x+1)^2dx}),
\end{eqnarray*}
$u_n$ denoting the solution of \eqref{1} emanating from $u_{n,0}$ at $t=0$.
\end{definition}

\begin{theorem}\label{global-pol}
For any $u_0 \in L^2_{(x+1)^2dx}$ and any $T>0$, there exists a
unique mild solution $u\in C([0,T];L^2_{(x+1)^2dx})\cap 
L^2(0,T;H^1_{(x+1)^2dx} )$ of \eqref{1}.
\end{theorem}

\noindent {\bf Proof.} We prove the existence and the uniqueness in two steps.\\
{\sc Step 1. Existence}\\
Since the embedding $L^2_b\subset L^2_{(x+1)^2 dx}$ is dense, for any given 
$u_0\in L^2_{(x+1)^2dx}$ we may construct 
a sequence $\{u_{n,0}\}\subset L^2_b$ such that $u_{n,0}\rightarrow
u_0$ in $L^2_{(x+1)^2 dx}$ as $n\rightarrow \infty$. For each $n$, let 
$u_n$ denote the solution of \eqref{1} 
emanating from $u_{n,0}$ at $t=0$, which is given by Theorem 
\ref{global-exp}. Then $u_n\in C([0,T];L^2_b)\cap L^2(0,T;H^1_b)$ and it 
solves 
\begin{eqnarray}
&& u_{n,t} + u_{n,x} + u_{n,xxx} + u_n u_{n,x} + a(x)u_{n} =0, \label{C1}\\
&& u_n(0,t) =0 \label{C2}\\
&& u_n(x,0)=u_{n,0}(x). \label{C3}
\end{eqnarray}
Multiplying \eqref{C1} by $(x+1)^2 u_n$ and integrating by parts, we obtain
\begin{eqnarray}
&&\frac{1}{2} \int_0^\infty (x+1)^2 |u_n (x,T)|^2 dx 
+3\int _0^T\!\!\!\int_0^\infty (x+1)|u_{n,x}|^2dxdt 
+\frac{1}{2}\int_0^T |u_{n,x}(0,t)|^2dt \nonumber\\
&&-\int_0^T\!\!\! \int_0^\infty (x+1)|u_n|^2 dxdt 
-\frac{2}{3}
\int_0^T\!\!\!\int_0^\infty (x+1)u_n^3 \, dxdt + 
\int_0^T\!\!\!\int_0^\infty (x+1)^2u_n^2a(x)dx \nonumber \\
&&\qquad = \frac{1}{2}\int_0^\infty (x+1)^2 |u_{n,0}(x)|^2 dx. \label{C4}
\end{eqnarray}
Scaling in \eqref{C1} by $u_n$ gives 
\begin{eqnarray*}
&&\frac{1}{2}\int_0^\infty |u_n(x,T)|^2 dx + 
\frac{1}{2}\int_0^T |u_{n,x}(0,t)|^2 dt + 
\int_0^T\!\!\!\int_0^\infty a(x)|u_n(x,t)|^2dxdt \\ 
&&=\frac{1}{2}\int_0^\infty |u_{n,0}(x)|^2dx,
\end{eqnarray*}
hence
\begin{equation}
||u_n||_{L^2(\R ^+)} \le ||u_{n,0}||_{L^2(\R ^+)} \le C
\label{C5}
\end{equation}
where $C=C( ||u_0||_{L^2(\R ^+)} )$.
It follows that 
\begin{eqnarray}
\frac{2}{3}\int_0^\infty (x+1)|u_n|^3dx 
&\le& \frac{2\sqrt{2}}{3} ||u_{n,x}||_{L^2(\R ^+)}^{\frac{1}{2}}
||u_n||_{L^2(\R ^+)}^{\frac{3}{2}} ||(x+1)u_n||_{L^2(\R ^+)} \nonumber \\
&\le& \int_0^\infty (x+1)|u_{n,x}|^2 dx + C \int_0^\infty (x+1)^2 |u_n|^2 dx
\label{C6}
\end{eqnarray}
which, combined to \eqref{C4}, gives
\begin{eqnarray*}
&&\frac{1}{2}\int_0^\infty (x+1)^2 |u_n(x,T)|^2 dx 
+ 2 \int_0^T\!\!\!\int_0^\infty (x+1)|u_{n,x}|^2 dxdt
+\frac{1}{2}\int_0^T |u_{n,x}(0,t)|^2 dt \\
&&\qquad \le \frac{1}{2} \int_0^\infty (x+1)^2 |u_{n,0}(x)|^2dx + 
C\int_0^T\!\!\!\int_0^\infty  (x+1)^2 |u_n (x,t)|^2 dxdt. 
\end{eqnarray*}
An application of Gronwall's lemma yields
\begin{eqnarray*}
||u_n||_{L^\infty (0,T;L^2_{(x+1)^2dx} )}
&\leq& C(T,||u_{n,0}||_{L^2_{ (x+1)^2dx}}),\\
%\label{C7}\\
||u_{n,x}||_{L^2(0,T;H^1_{ (x+1)^2 dx}   )}&\leq&
C(T,||u_{n,0}||_{L^2_{(x+1)^2 dx}}),\\
||u_{n,x}(0,.)||_{L^2(0,T)} &\leq& C(T,||u_{n,0}||_{L^2_{ (x+1)^2dx}}).
\end{eqnarray*}

Therefore, there exists a subsequence
of $\{u_n\}$, still denoted by $\{ u_n \}$, such that
\begin{equation}
\begin{cases}
u_n \rightharpoonup
u\,\,\mbox{weakly}\,\,\ast\,\,\mbox{in}\,\, L^\infty (0,T; L^2_{(x+1)^2dx}),\\
u_n \rightharpoonup
u\,\,\mbox{weakly}\,\,\mbox{in}\,\,
L^2(0,T; H^1_{(x+1)^2dx}),\\
u_{n,x}(0,.)  \rightharpoonup
u_x(0,.) \,\,\mbox{weakly}\,\,\mbox{in}\,\,
L^2(0,T).
\nonumber
\end{cases}
\end{equation}
Note that, for all $L>0$, $\{ u_n\}$ is bounded in 
$L^2(0,T;H^1(0,L))\cap H^1(0,T;H^{-2}(0,L))$, hence by Aubin's lemma, we have 
(after extracting a subsequence if needed)
$$u_n\to u\ \ \mbox{\rm strongly in } \  L^2(0,T;L^2(0,L)) \mbox{ for all } L>0.$$ 
This gives that $u_n u_{n,x} \to u u_x$ in the sense of distributions, hence the limit $u\in L^\infty (0,T;L^2_{(x+1)^2dx}) \cap L^2(0,T;H^1_{(x+1)^2dx})$ 
is a solution of \eqref{1}. 
Let us check that $u\in C([0,T]; L^2_{(x+1)^2dx})$. Since
$u\in C([0,T]; H^{-2}(\R ^+))\cap L^\infty (0,T; L^2_{(x+1)^2dx})$, we have that
$u\in C_w ([0,T]; L^2_{(x+1)^2dx})$ (see e.g. \cite{lions-magenes}), where
$C_w ([0,T]; L^2_{(x+1)^2 dx})$ denotes the space of sequentially weakly
continuous functions from $[0,T]$ into $L^2_{(x+1)^2dx}$.

We claim that $u\in L^3(0,T;L^3(\R ^+))$. Indeed, from Moser estimate (see \cite{Taylor3})
\begin{equation}
\label{Moser}
||u||_{L^\infty (\R ^+)}\le \sqrt{2}||u_x||^{\frac{1}{2}}_{L^2(\R ^+)}
||u||^{\frac{1}{2}}_{L^2(\R ^+)}
\end{equation}
and Young inequality we get
\begin{equation}
\label{exp-5bis}
\int_0^\infty |u|^3 dx 
\le ||u||_{L^\infty} ||u||^2_{L^2}
\le \sqrt{2} ||u_x||^{\frac{1}{2}}_{L^2} ||u||^{\frac{5}{2}}_{L^2}
\le \varepsilon ||u_x||_{L^2}^2 + c_{\varepsilon}
||u||_{L^2}^{\frac{10}{3}}
\end{equation}
where $\varepsilon >0$ is arbitrarily chosen and $c_\varepsilon$
denotes some positive constant.
Since 
$u\in C_w([0,T];L^2_{(x+1)^2dx}) \cap L^2(0,T;H^1_{(x+1)^2 dx})$, it follows that
$u\in L^3(0,T;L^3(\R ^+))$. On the other hand, $u(0,t)=0$ for $t\in (0,T)$ and
$u_x(0,.) \in L^2(0,T)$. 
Scaling in \eqref{1} by $(x+1)^2 u$ yields for all $t_1,t_2\in (0,T)$ 
\begin{eqnarray}
&&\frac{1}{2}\int_0^\infty (x+1)^2 |u(x,t_2)|^2 dx - 
\frac{1}{2}\int_0^\infty (x+1)^2 |u(x,t_1)|^2 dx \nonumber \\
&&=-3 \int_{t_1}^{t_2} \!\!\!\int_0^\infty (x+1)|u_x|^2 dxdt 
-\frac{1}{2}\int_{t_1}^{t_2}|u_x(0,t)|^2 dt +
\int_{t_1}^{t_2}\!\!\!\int_0^\infty (x+1)|u|^2 dxdt \nonumber\\
&&+\frac{2}{3}\int_{t_1}^{t_2}\!\!\!\int_0^\infty (x+1)u^3 dxdt 
-\int_{t_1}^{t_2}\!\!\!\int_0^\infty (x+1)^2 a(x) |u|^2 dxdt.
\label{LL28} 
\end{eqnarray}
Therefore $\lim_{t_1\to t_2} 
\left\vert ||u(t_2)||^2_{L^2_{(x+1)^2dx}}-||u(t_1)||^2_{L^2_{(x+1)^2dx}} 
\right\vert =0$.
Combined to the fact that $u\in C_w([0,T];L^2_{(x+1)^2dx})$, this yields
$u\in C([0,T],L^2_{(x+1)^2dx})$. \\
{\sc Step 2. Uniqueness}\\
Here, $C$ will denote a universal constant which may vary from line to line. 
Pick $u_0\in L^2_{(x+1)^2dx}$, and let $u,v\in C([0,T];L^2_{(x+1)^2dx})
\cap L^2(0,T;H^1_{(x+1)^2dx})$ be two mild solutions of \eqref{1}. Pick two
sequences $\{u_{n,0}\}$, $ \{v_{n,0} \}$ in $L^2_b$ for some $b>0$ such that 
\begin{eqnarray}
&&u_{n,0}\to u_0 \,\, \mbox{ strongly in }\,\,  L^2_{(x+1)^2dx}, \label{X1}\\
&&u_n \to u\,\,\mbox{weakly}\ast\,\,\mbox{in}\,\, 
L^\infty (0,T; L^2_{(x+1)^2dx}),\label{X2}\\
&&u_n \to u\,\,\mbox{weakly}\,\,\mbox{in}\,\,
L^2(0,T; H^1_{(x+1)^2dx})\label{X3}
\end{eqnarray}
and also
\begin{eqnarray}
&&v_{n,0}\to u_0 \,\, \mbox{ strongly in }\,\,  L^2_{(x+1)^2dx}, \label{X4}\\
&&v_n \to v\,\,\mbox{weakly}\ast\,\,\mbox{in}\,\, 
L^\infty (0,T; L^2_{(x+1)^2dx}),\label{X5}\\
&&v_n \to v\,\,\mbox{weakly}\,\,\mbox{in}\,\,
L^2(0,T; H^1_{(x+1)^2dx}).\label{X6}
\end{eqnarray}
We shall prove that $w=u-v$ vanishes on $\R ^+\times [0,T]$  by providing some 
estimate for $w_n=u_n-v_n$. Note first that $w_n$ solves the system
\begin{eqnarray}
&&w_{n,t} + w_{n,x}+ w_{n,xxx} + aw_n = f_n 
= v_n v_{n,x} - u_n u_{n,x}, \label{X7}\\
&&w_n(0,t)=0, \label{X8}\\
&&w_n(x,0)=w_{n,0}(x)=u_{n,0}(x)-v_{n,0}(x). \label{X9}
\end{eqnarray}
Scaling in \eqref{X7} by $(x+1)w_n$ yields 
\begin{eqnarray*}
&&\frac{1}{2}\int_0^\infty (x+1)|w_n(x,t)|^2dx + 
\frac{3}{2}\int_0^t \!\!\! \int_0^\infty |w_{n,x}|^2 dxd\tau 
-\frac{1}{2}\int_0^t\!\!\!\int_0^\infty |w_n|^2dxd\tau \nonumber\\
&&\le 
\frac{1}{2} \int_0^\infty (x+1)|w_{n,0}|^2 dx 
+ \int_0^t (\int _0^\infty (x+1)|w_n|^2dx)^{\frac{1}{2}}
 (\int _0^\infty (x+1)|f_n|^2dx)^{\frac{1}{2}}d\tau \\
&&\le \frac{1}{2} \int_0^\infty (x+1)|w_{n,0}|^2 dx 
+ \frac{1}{4} \sup_{0 < \tau < t}\int _0^\infty (x+1)|w_n(x,\tau )|^2dx \\
&&\qquad + [\int_0^T(\int _0^\infty (x+1)|f_n|^2dx)^{\frac{1}{2}} d\tau ]^2. 
\end{eqnarray*}
Since $||w_n(t)||_{L^2(\R ^+)}\le ||w_n(t)||_{L^2_{(x+1)dx}}$, this yields for $T<1/10$ 
\begin{eqnarray}
&&\sup_{0<t<T}\int_0^\infty (x+1)|w_n(x,t)|^2dx + \int_0^T\!\!\!\int_0^\infty 
|w_{n,x}|^2dxdt \nonumber\\
&&\qquad \le C[ \int_0^\infty (x+1)|w_{n,0}(x)|^2dx + 
\left( \int_0^T (\int_0^\infty (x+1) |f_n|^2 dx) ^{\frac{1}{2}}
d\tau \right)^2  ].
\label{X10}
\end{eqnarray}
It remains to estimate $\int_0^T(\int_0^\infty (x+1)|f_n|^2dx)^{\frac{1}{2}}dt$. 
We split $f_n$ into 
$$
f_n = (v_n-u_n) v_{n,x} + u_n ( v_{n,x} - u_{n,x} ) = f_n^1 + f_n^2.
$$
We have that 
\begin{eqnarray*}
\int_0^T (\int_0^\infty (x+1)|f_n^1|^2dx )^\frac{1}{2}dt 
&=& \int_0^T (\int_0^\infty (x+1) |w_n|^2 |v_{n,x}|^2dx )^{\frac{1}{2}} dt \\
&\le & \int_0^T ||w_n||_{L^\infty (\R ^+)}
(\int_0^\infty (x+1) |v_{n,x}|^2 dx)^{\frac{1}{2}} dt \\
&\le&  (\int_0^T ||w_n||^2_{L^\infty (\R ^+)}dt )^\frac{1}{2} (\int_0^T\!\!\!\int_0^\infty
(x+1) |v_{n,x}|^2dxdt )^\frac{1}{2}.
\end{eqnarray*}
By Sobolev embedding, we have that 
\begin{eqnarray*}
(\int_0^T ||w_n||^2_{L^\infty (\R ^+)} dt )^\frac{1}{2} 
&\le& 
(\int_0^T ||w_n||^2_{H^1 (\R ^+)}  dt  )^\frac{1}{2} \\
&\le& \sqrt{T} \sup_{0<t<T} ||w_n||_{L^2(\R ^+)} + 
||w_{n,x}||_{L^2(0,T;L^2(\R ^+))}\cdot
\end{eqnarray*}
Thus
\begin{eqnarray}
&&\int_0^T (\int_0^\infty (x+1)|f_n^1|^2dx)^{\frac{1}{2}} dt 
\le ||v_{n,x}||_{L^2(0,T;L^2_{(x+1)dx})}
\big( \sqrt{T} \sup_{0<t<T}||w_n||_{L^2(\R ^+)} \nonumber\\
&&\qquad+ ||w_{n,x}||_{L^2(0,T;L^2(\R ^+))}\big)
\label{X11}
\end{eqnarray}
On the other hand, we have that 
\begin{eqnarray}
&&\int_0^T ( \int_0^\infty (x+1)|f_n^2|^2 dx)^\frac{1}{2}dt \nonumber\\
&&\qquad = \int_0^T (\int_0^\infty (x+1)|u_n|^2 |w_{n,x}|^2 dx )^{\frac{1}{2}} dt\nonumber \\
&&\qquad \le \int_0^T ||(x+1)^{\frac{1}{2}} u_n||_{L^\infty (\R ^+)} 
||w_{n,x}||_{L^2(\R ^+)} dt \nonumber \\
&&\qquad \le C\int_0^T \big(||(x+1)^\frac{1}{2}u_n||_{L^2(\R ^+)}
+ ||(x+1)^\frac{1}{2} u_{n,x}||_{L^2( \R ^+)}
\big)||w_{n,x}||_{L^2(\R ^+ )} dt \nonumber \\
&&\qquad \le C \bigg(\sqrt{T} ||(x+1)u_n||_{L^\infty (0,T;L^2(\R ^+))}
\nonumber\\
&&\qquad\qquad +
||(x+1)^\frac{1}{2}u_{n,x}||_{L^2(0,T,L^2(\R ^+))}
\bigg)  ||w_{n,x}||_{L^2(0,T;L^2(\R ^+))}. \label{X12}
\end{eqnarray}
Gathering together \eqref{X10}, \eqref{X11} and \eqref{X12}, we conclude that 
for $T<1/10$
$$
h_n(T) \le K_n(T) h_n(T) + C||w_{n,0}||^2_{L^2_{(x+1)dx}} 
$$ 
where 
\begin{eqnarray}
h_n(t)&:=&\sup_{0<\tau <T}
\int_0^\infty (x+1)|w_n(x,\tau )|^2 dx
+\int_0^T \!\!\!\int_0^\infty |w_{n,x}|^2 dxdt\\
K_n(T)&\le& C \left( 
\int_0^T\!\!\!\int_0^\infty (x+1)|v_{n,x}|^2dxdt
+T ||(x+1)u_n||^2_{L^\infty (0,T;L^2(\R ^+))} \right. \nonumber \\ 
&&\left.  +\int_0^T\!\!\! \int_0^\infty (x+1)|u_{n,x}|^2 dxdt \right)
\end{eqnarray}
and $C$ denotes a universal constant. The following claim is needed.\\
{\sc Claim 3.}
$$\lim_{T\to 0}\limsup _{n\to \infty}\int_0^T\!\!\!\int_0^\infty
(x+1)|u_{n,x}|^2dxdt=0,\quad \lim_{T\to 0}\limsup _{n\to \infty}
\int_0^T\!\!\!\int_0^\infty (x+1)|v_{n,x}|^2dxdt=0.$$
Clearly, it is sufficient to prove the claim for the sequence 
$\{ u_n \}$ only. 
%Scaling in 
%$$
%u_{n,t} + u_{n,xxx} + u_{n,x} + u_n u_{n,x} + a(x)u_n=0 
%$$
%by $(x+1)u_n$ results in 
From \eqref{LL28} applied with $u=u_n$ on $[0,T]$, we obtain
\begin{eqnarray*}
&&\frac{1}{2}\int_0^\infty (x+1)^2|u_n(x,T)|^2dx +
%-\frac{1}{2}\int_0^\infty (x+1)|u_{n,0}|^2dx + 
3\int_0^T\!\!\!\int_0^\infty (x+1)|u_{n,x}|^2 dxdt  \nonumber \\
&&\le \frac{1}{2}\int_0^\infty (x+1)^2|u_{n,0}|^2dx  
+\int_0^T\!\!\!\int_0^\infty (x+1)|u_n|^2 dxdt 
+\frac{2}{3} \int_0^T\!\!\!\int_0^\infty  (x+1)|u_n|^3 dxdt.\label{P1}
\end{eqnarray*}
Combined to \eqref{C5}-\eqref{C6}, this gives
\begin{eqnarray}
&&||u_n(T)||^2_{L^2_{(x+1)^2dx}} +\int_0^T\!\!\!\int_0^\infty
(x+1)|u_{n,x}|^2dxdt \nonumber\\
&&\qquad \le ||u_{n,0}||^2_{L^2_{(x+1)^2 dx}} 
+C\int_0^T ||u_n||^2_{L^2_{(x+1)^2dx}}dt.
\label{P2}
\end{eqnarray}
It follows from Gronwall lemma that 
\begin{equation}
\label{gron}
||u_n(t)||^2_{L^2_{(x+1)^2dx}}\le ||u_{n,0}||^2_{L^2_{(x+1)^2dx}} e^{Ct}
\end{equation}
Using \eqref{gron} in \eqref{P2} 
and taking the limit sup as $n\to \infty$ gives for a.e. $T$
$$
||u(T)||^2_{L^2_{(x+1)^2dx}} + 
\limsup_{n\to \infty} \int_0^T\!\!\!\int_0^\infty |u_{n,x}|^2dxdt
\le e^{CT}||u_0||^2_{L^2_{(x+1)^2dx}} 
$$
As $u$ is continuous from $\R^+$ to $L^2_{(x+1)^2dx}$, we infer that
$$
\lim_{T\to 0} \limsup_{n\to \infty}\int_0^T\!\!\!\int_0^\infty 
|u_{n,x}|^2dxdt=0. 
$$
The claim is proved. Therefore, we have that
for $T>0$ small enough and $n$ large enough, $K_n(T)<\frac{1}{2}$, and hence
$$
h_n(T)\le 2C ||w_n(0)||^2_{L^2_{(x+1)dx}}.
$$
This yields
\begin{equation*}
||u-v||^2_{L^\infty (0,T;L^2_{(x+1)dx})} 
\le \liminf _{n\to \infty}
h_n(T) \le 2C \liminf _{n\to \infty} ||w_n(0)||^2_{L^2_{(x+1)dx}} =0
\end{equation*}
and $u=v$ for $0<t<T$. This proves the uniqueness for $T$ small enough.
The general case follows by a classical argument.
\qed

\begin{remark}
\begin{enumerate}
\item If we assume only that $u_0\in L^2_{(x+1)dx}$, then a proof similar to 
Step 1 gives the existence of a mild solution $u\in C([0,T];L^2_{(x+1)dx})
\cap L^2(0,T;H^1_{(x+1)dx})$ of \eqref{1}. The uniqueness of such a solution
is open. The existence and uniqueness of a solution issuing from
$u_0\in L^2_{(x+1)dx}$ in a class of functions involving a Bourgain
norm has been given in \cite{faminskii2}.
\item If $u_0\in L^2_{(x+1)^m dx}$ with $m\ge 3$, then 
$u\in C([0,T];L^2_{(x+1)^mdx})\cap L^2(0,T;H^1_{(x+1)^mdx})$ for all $T>0$
(see below Theorem \ref{dec-pol}). 
\end{enumerate}
\end{remark}
\section{Asymptotic Behavior}

\subsection{Decay in $L^2_{(x+1)^m dx}$}

\begin{theorem}
\label{dec-pol}
Assume that the function $a=a(x)$
satisfies (\ref{2}). Then, for all $R>0$ and $m\ge 1$, there exist
numbers $C > 0$ and $\nu > 0$ such that
$$||u(t)||_{L^2_{(x+1)^mdx}} \leq C\, e^{-\nu t}
||u_0||_{L^2_{(x+1)^mdx}} $$
for any solution given by Theorem \ref{global-pol}, whenever
$||u_0||_{ L^2_{(x+1)^mdx}  }\leq R$.
\end{theorem}

\noindent {\bf Proof.}
The proof will be done by induction in $m$.
We set
\begin{equation}\label{V0}
V_0(u) = E(u) = \frac{1}{2} \int_0^\infty u^2 dx
\end{equation}
and define the Lyapunov function $V_m$ for $m\ge 1 $ in an inductive way
\begin{equation}\label{pol1}
V_m(u) = \displaystyle\frac{1}{2}\int_0^\infty (x + 1)^m u^2 dx +
d_{m-1}V_{m-1}(u),
\end{equation}
where $d_{m-1} > 0$ is chosen sufficiently large (see below).

Suppose first that $m=1$ and put $V = V_1$. Multiplying the first equation in (\ref{1})
by $u$ and integrating by parts over $\R^+ \times (0,T)$, we obtain
\begin{equation}\label{exp-3}
\frac{1}{2}\int_0^\infty |u(x,T)|^2 dx = \frac{1}{2}\int_0^\infty
|u_0(x)|^2 dx - \int_0^T\int_0^\infty a(x)|u|^2 dxdt
- \frac{1}{2}\int_0^T u_x^2(0,t)dt.
\end{equation}
Now, multiplying the equation by $xu$, we deduce that
\begin{eqnarray}
&&\vspace{1mm}\displaystyle\frac{1}{2}\int_0^\infty x |u(x,T)|^2 dx -
\frac{1}{2}\int_0^\infty x |u_0(x)|^2 dx + \frac{3}{2}
\int_0^T\int_0^\infty u_x^2 dx dt \nonumber \\
&& -\displaystyle\frac{1}{2}
\int_0^T\int_0^\infty u^2 dx dt-\displaystyle\frac{1}{3}
\int_0^T\int_0^\infty u^3 dx dt+\int_0^T\int_0^\infty x a(x)|u|^2dxdt
= 0.\qquad
\label{exp-4}
\end{eqnarray}
Combining (\ref{exp-3}) and (\ref{exp-4}) it follows that
\begin{eqnarray}
&& \vspace{1mm} V(u) - V(u_0) + (d_0+1)\left( \displaystyle\frac{1}{2}\int_0^T
u_x^2(0,t)dt + \int_0^T\int_0^\infty a(x)|u|^2 dxdt\right) \nonumber\\
&& \vspace{1mm} + \displaystyle \frac{3}{2}
\int_0^T\int_0^\infty u_x^2 dx dt-\displaystyle\frac{1}{2}
\int_0^T\int_0^\infty u^2 dx dt -\displaystyle\frac{1}{3}
\int_0^T\int_0^\infty u^3 dx dt \nonumber\\
&& \qquad\qquad\qquad\qquad\qquad\qquad\qquad\qquad + \displaystyle\int_0^T\int_0^\infty xa(x)|u|^2 dxdt = 0.
\label{exp-5}
\end{eqnarray}

The next step is devoted to estimate the nonlinear term in the left
hand side of (\ref{exp-5}). To do that, we first assume that
$||u_0||_{L^2} \le 1$.

By \eqref{exp-5bis} we have that 
\begin{equation*}
\int_0^\infty |u|^3 dx 
\le \varepsilon ||u_x||_{L^2}^2 + c_{\varepsilon}
||u||_{L^2}^{\frac{10}{3}}
\end{equation*}
for any $\varepsilon >0$  and some constant $c_\varepsilon >0$.
Thus, if $||u_0||_{L^2} \le 1$, we have
$||u||_{L^2}^{\frac{10}{3}} \leq ||u||_{L^2}^2$ and
\begin{equation}
\label{exp-6}
\int_0^T\int_0^\infty |u|^3 dx dt \leq \varepsilon
\int_0^T\int_0^\infty u_x^2 dx dt+
c_\varepsilon\int_0^T\int_0^\infty u^2 dx dt.
\end{equation}
Moreover, according to \cite{linares-pazoto}, there exists $c_1>0$,
satisfying
\begin{equation}\label{exp-7}
\int_0^T\int_0^\infty u^2 dx dt \leq c_1\{ \frac{1}{2}
\int_0^T u_x^2(0,t) dt+
\int_0^T\int_0^\infty a(x) u^2 dx dt \}.
\end{equation}
Now, combining (\ref{exp-5})-(\ref{exp-7}) and taking $\varepsilon <
\frac{1}{2}$ and $d_0 := 2 c_1(\frac{1}{2} +
\frac{c_\varepsilon}{3})$ we obtain
\begin{eqnarray}
&&V(u(T)) - V(u_0) +
\frac{d_0+1}{2}(\frac{1}{2}\int_0^T u_x^2(0,t)dt +
\int_0^T\int_0^\infty a(x)|u|^2 dxdt)
\nonumber\\
&&\qquad+\,(\frac{3}{2} - \frac{\varepsilon}{3})
\int_0^T\int_0^\infty u_x^2 dx dt
+\displaystyle\int_0^T\int_0^\infty xa(x)|u|^2 dxdt
%+\frac{1}{2}\int_0^T u_x^2(0,t)dt
\leq 0
\label{exp-71}
\end{eqnarray}
or
\begin{equation}
\label{exp-72}
V(u(T)) - V(u_0) \leq -
\widetilde{c}\,\{\displaystyle\int_0^T u_x^2(0,t)dt +
\int_0^T\int_0^\infty (x + 1) a(x)|u|^2 dxdt + \int_0^T\int_0^\infty
u_x^2 dx dt \}
\end{equation}
where $\widetilde{c}>0$. We aim to prove the existence  of 
a constant $c>0$ satisfying
\begin{equation}
\label{exp-73}
V(u(T)) - V(u_0) \leq - c\,V(u_0)
\end{equation}
Indeed, such an inequality gives at once the decay $V(u(t))\le c e^{-\nu t} 
V(u_0)$. To this end, we need to establish two claims. 

{\sc Claim 4.} There exists $c>0$ such that
$$\displaystyle\int_0^T V(u) dt \leq c\, \{\int_0^T u_x^2(0,t) dt
+ \int_0^T\int_0^\infty (x +1) a(x) u^2 dx dt\}.$$

Since $u_0\in L^2_{(x+1)dx}\subset L^2$, from (\ref{2}) and
(\ref{exp-7}) we get
\begin{eqnarray*}
\int_0^T V(u) dt
&=&
\frac{1}{2}\int_0^T\int_0^\infty (x + 1) u^2 dx
dt + \frac{d_0}{2} \int_0^T\int_0^\infty u^2 dx dt \\
&\leq& \frac{c_1 d_0}{2}
\{\displaystyle \frac{1}{2} \int_0^T
u_x^2(0,t) dt + \int_0^T\int_0^\infty a(x) u^2 dx dt\}\\
&&\qquad + \frac{1}{2}\int_0^T\int_0^{x_0} (x + 1 )u^2 dx dt
+\frac{1}{2}\int_0^T\int_{x_0}^\infty (x + 1) u^2 dx dt\\
&\leq&  \frac{c_1 d_0}{2}
\{\frac{1}{2} \int_0^T u_x^2(0,t) dt +
\int_0^T\int_0^\infty a(x) u^2 dx dt\} \\
&&\qquad  +\frac{1}{2}(x_0 +
1)\int_0^T\int_0^{x_0}u^2 dx dt +
\frac{1}{2}\int_0^T\int_{x_0}^\infty (x + 1)\frac{a(x)}{a_0} u^2 dx dt\\
&\leq& c\,\{ \int_0^T u_x^2(0,t) dt +
\int_0^T\int_0^\infty (x + 1) a(x) u^2 dx dt\}.
\end{eqnarray*}

{\sc Claim 5.}
\begin{equation}
\label{exp-74}
V(u_0) \leq C (\int_0^T u_x^2(0,t) dt
+ \int_0^T\int_0^\infty (x + 1) a(x) u^2 dx dt
+ \int_0^T\int_0^\infty u_x^2 dx dt )
\end{equation}
where $C > 0$.

Multiplying the first equation in (\ref{1}) by $(T-t)u$ and
integrating by parts in $(0,\infty)\times (0,T)$, we obtain
\begin{equation}\label{exp-75}
\begin{array}{l}
\vspace{1mm}\displaystyle\frac{T}{2}\int_0^\infty |u_0(x)|^2 dx =\\
\displaystyle\frac{1}{2}\int_0^T\int_0^\infty |u|^2 dx dt + \int_0^T\int_0^\infty
(T-t) a(x)|u|^2 dxdt +\frac{1}{2}\int_0^T (T-t) u_x^2(0,t)dt,
\end{array}
\end{equation}
and therefore, using (\ref{exp-7})
\begin{equation}\label{exp-8}
\begin{array}{l}
\displaystyle\int_0^\infty |u_0(x)|^2 dx
\leq C\left( \int_0^T\int_0^\infty a(x)|u|^2 dxdt
+ \int_0^T u_x^2(0,t) dt\right) .
\end{array}
\end{equation}
Now, multiplying by $(T-t)x u$, it follows that
\begin{equation}
\begin{array}{l}
\vspace{1mm}-\displaystyle\frac{T}{2}\int_0^\infty x |u_0(x)|^2 dx +
\frac{1}{2}\int_0^T\int_0^\infty x |u|^2 dx dt
+\,\frac{3}{2}\displaystyle\int_0^T\int_0^\infty (T-t) u^2_x dx dt\\
\vspace{1mm}- \displaystyle\frac{1}{2}\int_0^T\int_0^\infty (T- t)
u^2 dx dt +\displaystyle\int_0^T\int_0^\infty (T - t)x a(x)|u|^2 dx
dt - \\
\qquad\qquad\qquad\qquad\qquad\qquad\qquad\qquad\qquad
-\displaystyle\frac{1}{3}\int_0^T \int_0^\infty (T - t) u^3 dx dt= 0.
\nonumber
\end{array}
\end{equation}
The identity above and (\ref{exp-6}) allow us to conclude that
\begin{equation}
\begin{array}{l}\label{exp-9}
\vspace{1mm}\displaystyle\int_0^\infty x |u_0(x)|^2 dx \\
\vspace{1mm}\leq C\,\{ \displaystyle\int_0^T\!\!\!\int_0^\infty (x + 1)
|u|^2 dx dt +\,\displaystyle\int_0^T\!\!\!\int_0^\infty u^2_x dx
dt+\displaystyle\int_0^T\!\!\!\int_0^\infty x a(x)|u|^2 dx dt +\\
+\displaystyle\int_0^T\!\!\! \int_0^\infty  |u|^3 dx dt\}
\leq C\,\{ \displaystyle\int_0^T V(u(t))dt
+ \int_0^T\!\!\!\int_0^\infty x a(x) u^2 dx dt
+\,\displaystyle\int_0^T\!\!\!\int_0^\infty u^2_x dx dt\}
\end{array}
\end{equation}
for some $C > 0$. Claim 5 follows
 from Claim 4 and (\ref{exp-8})-(\ref{exp-9}). \qed
%On the other hand
%\begin{eqnarray}
%\int_0^T \int_0^{\infty} (x+1)|u|^2 dxdt
%&\leq & (x_0+1)\int_0^T\int_0^{x_0} |u|^2 dxdt
%+\int_0^T\int_{x_0}^\infty \frac{a(x)}{a_0} (x+1)|u|^2 dxdt
%\nonumber\\
%&\leq & C
%\left(
%\int_0^T |u_x(0,t)|^2 dt + \int_0^T\int_0^\infty a(x) (x+1)|u|^2 dxdt
%\right)
%\label{exp-9bis}
%\end{eqnarray}

The previous computations give us \eqref{exp-73} (and the exponential decay) 
when \\
$||u_0||_{L^2}
\leq 1$. The general case is proved as follows. Let $u_0\in
L^2_{(x+1)dx} \subset L^2$ be such that $||u_0||_{L^2} \leq R$. Since $u\in
C(\mathbb{R} ^+ ;L^2(\mathbb{R} ^+))$ and $||u(t)||_{L^2} \leq \alpha
e^{-\beta t}||u_0||_{L^2}$, where $\alpha =\alpha (R)$ and
$\beta =\beta (R)$ are
positive constants, $||u(T)||_{L^2} \leq 1$ if we pick $T$
satisfying $\alpha e^{-\beta T} R< 1$. Then, it follows from
(\ref{exp-5})-(\ref{exp-5bis}) and \eqref{exp-73}  that
for some constants $\nu >0, \ c>0,\  C>0$
$$V(u(t+T)) \leq c e^{-\nu t}V(u(T)) \leq c (T||u_0||^2_{L^2} + T||u_0||_{L^2}^{\frac{10}{3}} + V(u_0)) e^{-\nu t},$$
hence
$$V(u(t)) \leq C e^{-\nu t}V(u_0),$$
where $C =  C(R)$, which concludes the proof when $m=1$.

\vglue 0.2 cm

\noindent{\bf Induction Hypothesis:} There exist $c > 0$ and  $\rho >0$ such that
if $V_{m-1}(u_0)\le \rho$, we have
\begin{equation}
\begin{array}{l}
\vspace{1mm}V_{m}(u) - V_{m}(u_0) \hfill (*)_{m}\\
\vspace{1mm}\leq - c\{\displaystyle\int_0^T u_x^2(0,t)
dt + \int_0^T\int_0^\infty (x + 1)^{m-1}u_x^2 dx dt +
\int_0^T\int_0^\infty (x + 1)^{m} a(x) u^2dx dt \}\\
V_{m}(u_0) \hfill (**)_{m}\\
\leq c\, \{\displaystyle\int_0^T u_x^2(0,t) dt +
\int_0^T\int_0^\infty (x + 1)^{m-1}u_x^2 dx dt +
\int_0^T\int_0^\infty (x + 1)^{m} a(x) u^2 dx dt \}.\nonumber
\end{array}
\end{equation}
By \eqref{exp-72}-\eqref{exp-74}, the induction hypothesis is true for
$m=1$. Pick now an index $m\ge 2$ and assume that $d_0, ..., d_{m-2}$ have been constructed so that
$(*)_k - (**)_k$ are fulfilled for $1\le k \le m-1$. We aim to prove that
for a convenient choice of the constant $d_{m-1}$ in \eqref{pol1}, the
properties $(*)_m-(**)_m$ hold true.

Let us investigate first $(*)_m$.
We multiply the first equation in \eqref{1}  by $(x + 1)^m u$ to obtain
\begin{equation}
\begin{array}{l}\label{exp-10}
\vspace{2mm}V_m(u) - V_m(u_0)  - d_{m-1}(V_{m-1}(u) - V_{m-1}(u_0))\\
-\displaystyle\frac{m(m-1)(m-2)}{2}\int_0^T\int_0^\infty (x +
1)^{m-3} u^2 dx dt + \displaystyle\frac{1}{2}\displaystyle\int_0^T u_x^2(0,t)
dt\\
\vspace{2mm}+\displaystyle\frac{3m}{2}\int_0^T\int_0^\infty (x +
1)^{m-1} u^2_x dx dt - \displaystyle\frac{m}{2}\int_0^T\int_0^\infty (x + 1)^{m-1}u^2 dx dt\\ -
\displaystyle\frac{m}{3}\int_0^T\int_0^\infty (x + 1)^{m-1}u^3 dx dt
+ \displaystyle\int_0^T\int_0^\infty (x + 1)^m a(x) u^2 dx dt = 0.
\end{array}
\end{equation}

The next steps are devoted to estimate the terms in the above
identity. First, combining (\ref{2}) and (\ref{exp-7}) we infer the
existence of a positive constant $c > 0$ such that
\begin{equation}\label{exp-12}
\begin{array}{l}
\vspace{1mm}\displaystyle\int_0^T\int_0^\infty (x + 1)^{m-1} u^2 dx dt \\
\vspace{1mm}= \displaystyle\int_0^T\int_0^{x_0} (x + 1)^{m-1} u^2 dx dt + \int_0^T\int_{x_0}^\infty
(x + 1)^{m-1} u^2 dx dt \\
\vspace{1mm}\leq (x_0 + 1)^{m-1}\displaystyle\int_0^T\int_0^\infty u^2 dx dt +
\displaystyle\frac{1}{a_0}\int_0^T\int_0^\infty a(x) (x + 1)^{m-1} u^2
dx dt \\
\vspace{1mm}\leq c\,\{\displaystyle\int_0^T u_x^2 (0,t) dt + \displaystyle\int_0^T\int_0^\infty (x + 1)^{m-1} a(x)
u^2 dx dt\} \\
\leq - c\,\{V_{m-1}(u) - V_{m-1}(u_0)\}
\end{array}
\end{equation}
where we used $(*)_{m-1}$.
In the same way
\begin{equation}
\begin{array}{l}\label{exp-13}
\vspace{1mm}\displaystyle\int_0^T\int_0^\infty (x + 1)^{m-3} u^2 dx dt \\
\leq \displaystyle\int_0^T\int_0^\infty (x + 1)^{m-1} u^2 dx dt
\leq - c\,\{V_{m-1}(u) - V_{m-1}(u_0)\}
\end{array}
\end{equation}
where $c > 0$ is a positive constant. Moreover, assuming
$V_{m-1}(u_0) \le \rho$ with $\rho >0$ small enough (so that by exponential 
decay of $V_{m-1}(u(t))$ we have $\int_0^\infty (x+1)^{m-1}|u(x,t)|^2dx \le 1$
for all $t\ge 0$) 
and proceeding as in the case $m=1$, we obtain the existence of
$\varepsilon > 0$ and $c_\varepsilon > 0$ satisfying
\begin{equation}\label{exp-14}
\begin{array}{l}
\vspace{1mm}\displaystyle\int_0^T\int_0^\infty (x + 1)^{m-1} |u|^3 dx
dt \\
\leq \varepsilon\displaystyle\int_0^T\int_0^\infty (x + 1)^{m-1}
u^2_x dx dt + c_{\varepsilon}\int_0^T\int_0^\infty (x + 1)^{m-1} u^2
dx dt.
\end{array}
\end{equation}
Indeed,
\begin{eqnarray}
&&\displaystyle\int_0^\infty (x + 1)^{m-1} |u|^3 dx \\
&&\leq ||u||_{L^\infty}\int_0^\infty (x + 1)^{m-1} u^2 dx
\leq \sqrt{2} ||u_x||_{L^2}^{\frac{1}{2}}||u||_{L^2}^{\frac{1}{2}}
\int_0^\infty (x+1)^{m-1}u^2 dx \nonumber \\
&&\leq \varepsilon\displaystyle\int_0^\infty (x + 1)^{m-1} u^2_x dx
+ c_\varepsilon \int_0^\infty u^2dx
+ c_{\varepsilon}\left(\int_0^\infty (x + 1)^{m-1} u^2 dx
\right)^{2}.\nonumber
\end{eqnarray}
Then, if we return to (\ref{exp-10}) and take $\varepsilon <9/2$
and $d_{m-1} > 0$ large enough,
from (\ref{exp-12})-(\ref{exp-14}) if follows that
\begin{equation}
\begin{array}{l}\label{exp-15}
\vspace{2mm}V_m(u) - V_m(u_0) \\
\vspace{2mm}\leq
-c\,\{\displaystyle
\int_0^T u^2_x (0,t) dt +
\int_0^T\int_0^\infty (x + 1)^{m-1} u^2_x dx dt +
\displaystyle
 \int_0^T\int_0^\infty a(x) (x + 1)^m u^2 dx dt\}\\
+\,
\displaystyle \frac{d_{m-1}}{2} (V_{m-1}(u) - V_{m-1}(u_0)).
\end{array}
\end{equation}
This yields $(*)_m$, by $(*)_{m-1}$. Let us now check $(**)_m$.
It remains to estimate the terms in the right
hand side of (\ref{exp-15}). We multiply the first equation
in \eqref{1} by $(T-t)(x + 1)^m u$ to obtain
\begin{equation}
\begin{array}{l}
\vspace{1mm}
\displaystyle
\frac{T}{2}\int_0^\infty (x + 1)^m u_0^2 dx = \frac{1}{2}\int_0^T\int_0^\infty (x + 1)^m
u^2 dx dt\\
\vspace{1mm}-\displaystyle\frac{m(m-1)(m-2)}{2}\int_0^T\int_0^\infty 
(T-t)(x + 1)^{m - 3}u^2 dx dt + \displaystyle\frac{1}{2}\int_0^T
(T-t)u_x^2(0,t) dt
 \\
\vspace{1mm}+\displaystyle\frac{3m}{2}\int_0^T\int_0^\infty (T - t)
(x + 1)^{m-1} u^2_x dx dt - \frac{m}{2}\int_0^T\int_0^\infty (T - t) (x +
1)^{m-1}u^2 dx dt\\
- \displaystyle\frac{m}{3}\int_0^T\int_0^\infty (T - t) (x +
1)^{m-1}u^3 dx dt + \int_0^T\int_0^\infty (T - t)(x + 1)^m a(x) u^2 dx dt. \nonumber
\end{array}
\end{equation}
Then, proceeding as above, we deduce that
\begin{equation}
\begin{array}{l}
\vspace{1mm}\displaystyle\int_0^T (x + 1)^m u_0^2 dx \\
\vspace{1mm}\leq c\,\{\displaystyle\int_0^T\int_0^\infty (x + 1)^{m-1} u^2 dx dt + \int_0^T
u_x^2(0,t) dt + \int_0^T\int_0^\infty (x + 1)^{m-1}u_x^2 dx dt\\
\vspace{1mm}+ \displaystyle\int_0^T\int_0^\infty (x + 1)^m a(x) u^2 dx dt\}\\
\leq c\{\displaystyle\int_0^T u_x^2(0,t) dt +
\int_0^T\int_0^\infty (x + 1)^{m-1}u_x ^2 dx dt +
\displaystyle\int_0^T\int_0^\infty (x + 1)^m a(x) u^2 dx dt
\}.\nonumber
\end{array}
\end{equation}
Combined to $(**)_{m-1}$, this yields $(**)_m$. This completes the construction
of the sequence $\{ V_m\}_{m\ge 1}$ by induction.

Let us now check the exponential decay of $V_m$ for $m\ge 2$. It follows from
$(*)_m-(**)_{m}$ that
$$V_m(u) - V_m(u_0) \leq - c\,V_m(u_0)$$
where $c > 0$, which completes the proof when $V_{m-1}(u_0)\le \rho$. The global result ($V_{m-1} (u_0)\leq R$) is obtained as above for $m=1$.
\qed

\begin{corollary}\label{c3} 
Let $a=a(x)$ fulfilling (\ref{2}) and $a\in W^{2,\infty }(0,\infty)$. Then for any 
$R>0$, there exist positive constants $c=c(R)$ and $\mu = \mu (R)$ such that 
%and, in addition, $\frac{d^2 a}{d x^2}\in L^\infty$
\begin{equation}
||u_x(t)||_{L^2(\R ^+)} \leq c \frac{e^{-\mu t}}{\sqrt{t}}
||u_0||_{L^2_{(x+1)dx}}
\label{L4}
\end{equation}
for all $t>0$ and all $u_0\in L^2_{(x+1)dx}$ satisfying 
$|| u_0 ||_{L^2_{(x+1)dx}}\le R$. 
\end{corollary}

\noindent {\bf Proof.}
Pick any $R>0$ and any $u_0\in L^2_{(x+1)dx}$ with 
$||u_0||_{L^2_{(x+1)dx}}\le R$. By Theorem \ref{dec-pol}
there are some constants $C=C(R)$ and $\nu =\nu (R)$ such that
\begin{equation}
\label{L4bis}
|| u(t) ||_{L^2_{(x+1)dx}} \le C e^{-\nu t} ||u_0||_{L^2_{(x+1)dx}}.
\end{equation}
Using the multiplier $t(u^2 + 2 u_{xx})$ we obtain after some integrations by parts 
that for all $0<t_1<t_2$ 
\begin{eqnarray}
&&t_2\int_0^\infty u_x^2(x,t_2) dx +
\int_{t_1}^{t_2} tu_x^2 (0,t) dt + 2\int_{t_1}^{t_2}\!\!\!\int_0^\infty t
a(x)u_x^2 dx dt + \int_{t_1}^{t_2} t u_{xx}^2 (0,t) dt\nonumber\\
&&=-\frac{1}{3}\int_{t_1}^{t_2}\!\!\!\int_0^\infty u^3 dx
dt + \frac{t_2}{3}\int_0^\infty u^3(x,t_2) dx + 
\int_{t_1}^{t_2}\!\!\!\int_0^\infty t u^3
a(x)dx dt \nonumber\\
&&+\int_{t_1}^{t_2}\!\!\!\int_0^\infty u_x^2 dxdt +
\int_{t_1}^{t_2}\!\!\!\int_0^\infty t a''(x) u^2 dx dt.
\label{L5}
%\leq
%\varepsilon(1+T ||a||_{L^\infty})\displaystyle\int_0^T\int_0^\infty u_x^2 dx dt +
%c_\varepsilon (1 + T ||a||_{L^\infty })
%\displaystyle \int_0^T ||u||^{\frac{10}{3}}_{L^2} dt  +
%\frac{Tc_\varepsilon}{3}||u(T)||^{\frac{10}{3}}_{L^2} \\+
%\displaystyle\frac{T\varepsilon}{3}\int_0^\infty u_x^2(x,T) dx
%+\displaystyle \int_0^T\int_0^\infty u_x^2 dx dt
%+ c T ||a''||_{L^\infty}\displaystyle \int_0^T\int_0^\infty u^2 dx dt
%\end{array}
\end{eqnarray}
1. Let us assume first that $T>1$. Applying \eqref{L5} on the time interval
$[T-1,T]$, we infer that
\begin{equation}
\int_0^\infty |u_x(x,T)|^2 dx
\leq c\left(\int_{T-1}^T\!\int_0^\infty |u|^3 dx dt + 
||u(T)||^3_{L^3(\R ^+)} + \int_{T-1}^T||u||^2_{H^1(\R ^+)} dt\right) .
\label{L8}
\end{equation}
To estimate the cubic terms in \eqref{L8}, we use \eqref{exp-5bis} to obtain
\begin{eqnarray}
&&\int_0^\infty |u_x(x,T)|^2dx \le \varepsilon \int_0^\infty |u_x(x,T)|^2dx 
\nonumber\\
&&\qquad + c_\varepsilon\big( ||u(T)||_{L^2(\R ^+)}^{\frac{10}{3}} + \int_{T-1}^T 
(||u||^2_{H^1(\R ^+)} + ||u||_{L^2(\R ^+)}^{\frac{10}{3}})dt\big).
\label{L9} 
\end{eqnarray}
Note that by \eqref{L4bis}
$$
||u(T)||_{L^2(\R ^+)}^{\frac{10}{3}} 
\le (C e^{-\nu T}||u_0||_{L^2_{(x+1)dx}} )^{\frac{10}{3}}
\le C^{\frac{10}{3}} R^{\frac{4}{3}} e^{-\nu T}||u_0||^2_{L^2_{(x+1)dx}} .$$
It follows from \eqref{exp-5}, \eqref{exp-5bis}, and \eqref{L4bis} that 
\begin{eqnarray}
&&\int_{T-1}^T(||u||^2_{H^1(\R ^+)} + ||u||^{\frac{10}{3}}_{L^2(\R ^+)})dt 
\nonumber\\
&&\qquad \le C\left( V_1(u(T-1))+ \int_{T-1}^T
\big( ||u||^2_{L^2 (\R ^+)} +||u||^{\frac{10}{3}}_{L^2(\R ^+)}\big) dt \right)
\nonumber\\
&&\qquad \le C e^{-\nu T}||u_0||^2_{L^2_{(x+1)dx}}\label{L10}
\end{eqnarray}
where $C=C(R,\nu)$. 
\eqref{L4} for $T\ge 1$ follows from \eqref{L9} and \eqref{L10} by choosing
$\varepsilon <1$ and $\mu < \nu$.\\
2. Assume now that $T\le 1$.  Estimating again the cubic terms in \eqref{L5} 
(with $[t_1,t_2]=[0,T]$) by using \eqref{exp-5bis}, we obtain 
\begin{eqnarray}
T\int_0^\infty u_x^2(x,T)dx 
&\le& \frac{T}{3}
\left( \varepsilon ||u_x(T)||^2_{L^2(\R ^+ )}
+C_\varepsilon ||u(T)||^{\frac{10}{3}}_{L^2(\R ^+)} \right) \nonumber\\
&&\quad +C_\varepsilon \int_0^T (||u||^2_{H^1(\R ^+)}
+||u||^{\frac{10}{3}}_{L^2(\R ^+)}) dt.
\label{L6}
\end{eqnarray}
By \eqref{exp-5}, \eqref{exp-5bis} and \eqref{L4bis}, we have that
\begin{equation}
\int_0^1\!\!\!\int_0^\infty |u_x|^2dxdt \le C(R) ||u_0||^2_{L^2_{(x+1)dx}}
\label{L7}
\end{equation}
which, combined to \eqref{L6} with $\varepsilon =1$ and \eqref{L4bis}, gives 
$$
||u_x(T)||^2_{L^2(\R ^+)} \le C(R)T^{-1}||u_0||^2_{L^2_{(x+1)dx}}
$$
for all $T<1$. This gives \eqref{L4} for $T<1$. 
\qed
Corollary \ref{c3} may be extended (locally) to the weighted space $L^2_{(x+1)^mdx}$ ($m\ge 2$) in following the method of proof 
of \cite[Theorem 1.1]{pazoto-rosier}.

\begin{corollary}\label{c4} 
Let $a=a(x)$ fulfilling \eqref{2} and $m\ge 2$. Then there exist some
constants $\rho >0$, $C>0$ and $\mu >0$ such that  
%and, in addition, $\frac{d^2 a}{d x^2}\in L^\infty$
$$||u(t)||_{H^1_{(x+1)^mdx}} \leq C \frac{e^{-\mu t}}{\sqrt{t}}
||u_0||_{L^2_{(x+1)^mdx}}$$
for all $t>0$ and all $u_0\in L^2_{(x+1)^mdx}$ satisfying 
$|| u_0 ||_{L^2_{(x+1)^mdx}}\le \rho$. 
\end{corollary}
\noindent {\bf Proof.}
We first prove estimates for the linearized problem
\begin{eqnarray}
&&u_t+u_x+u_{xxx}+au=0 \label{W1}\\
&&u(0,t)=0 \label{W2} \\
&&u(x,0)=u_0(x) \label{W3}
\end{eqnarray}
and next apply a perturbation argument to extend them to the nonlinear
problem \eqref{1}. Let us denote by $W(t)u_0=u(t)$ the solution of \eqref{W1}-\eqref{W3}. 
By computations similar to those performed in the proof of Theorem \ref{dec-pol}, we have that 
\begin{equation*}
||W(t)u_0||_{L^2_{(x+1)^mdx}} \le C_0 e^{-\nu t}||u_0||_{L^2_{(x+1)^mdx}}.
%\int_0^T ||W(t)u_0||^2_{H^1_{(x+1)^mdx}} dt \le C_2||u_0||^2_{L^2_{(x+1)^m dx}}.
\end{equation*}
We need the\\ 
{\sc Claim 6.} Let $k\in \{ 0, ... , 3\}$. Then there exists a constant $C_k>0$ 
such that for any $u_0\in H^k_{ (x+1)^m dx}$, 
\begin{equation}
||W(t)u_0||_{H^k_{(x+1)^m dx}} \le C_k e^{-\nu t} 
||u_0||_{H^k_{(x+1)^m dx}}.
\label{L12} 
\end{equation}
Indeed, if $u_0\in H^3_{(x+1)^mdx}$, then $u_t(.,0)\in L^2_{(x+1)^{m-3}dx}$, and
since $v=u_t$ solves \eqref{W1}-\eqref{W2}, we also have that 
$$
||u_t(.,t)||_{L^2_{ (x+1)^{m-3} dx}} \le C_0 e^{-\nu t} 
||u_t(.,0)||_{L^2_{(x+1)^{m-3}dx}}.
$$
Using \eqref{W1}, this gives 
$$
||W(t)u_0||_{H^3_{(x+1)^m dx}} \le C_3 e^{-\nu t} 
||u_0||_{H^3_{(x+1)^mdx}}.
$$
This proves \eqref{L12} for $k=3$. The fact that \eqref{L12} is valid for
$k=1,2$ follows from a standard interpolation argument, for 
$H^k_{(x+1)^mdx}=[H^0_{(x+1)^m dx},H^3_{(x+1)^mdx}]_{\frac{k}{3}}$.

\begin{lemma}\label{l2}
Pick any number $\mu \in (0,\nu)$. Then there exists some constant
$C=C(\mu )>0$ such that for any $u_0\in L^2_{(x+1)^mdx}$ 
\begin{equation}
||W(t)u_0||_{H^1_{(x+1)^mdx}} \le C \frac{e^{-\mu t}}{\sqrt{t}}
||u_0||_{L^2_{(x+1)^mdx}}\cdot
\label{W100}
\end{equation}
\end{lemma}
\noindent {\bf Proof.}
Let $u_0\in L^2_{(x+1)^m dx}$ and set $u(t)=W(t)u_0$ for all $t\ge 0$. 
By scaling in \eqref{W1} by $(x+1)^mu$, we see that for some constant
$C_K=C_K(T)$ 
$$
||u||_{L^2(0,1;H^1_{(x+1)^mdx})} \le C_K 
||u_0||_{L^2_{(x+1)^m dx}}\cdot
$$
This implies that $u(t)\in H^1_{(x+1)^mdx}$ for a.e. $t\in (0,1)$ which, 
combined to \eqref{L12}, gives that $u(t)\in H^1_{(x+1)^m dx}$ for all $t>0$.
Pick any $T\in (0,1]$. Note that, by \eqref{L12},  
\begin{equation}\label{W4}
||u(T)||_{H^1_{(x+1)^mdx}} \leq C_1 e^{-\nu (T-t)} 
||u(t)||_{H^1_{(x+1)^mdx}}, \quad \forall t\in (0,T).
\end{equation}
Integrating with respect to $t$ in \eqref{W4} yields
$$[C_1^{-1}||u(T)||_{H^1_{(x+1)^mdx}}]^2\int_0^T e^{2\nu (T-t)}dt 
\leq \int_0^T ||u(t)||^2_{H^1_{(x+1)^mdx}}dt,$$
and hence
\begin{eqnarray*}
||u(T)||_{H^1_{(x+1)^mdx}} &\leq& 
C_K\,C_1 \sqrt{\frac{2\nu}{e^{2\nu T}-1}}
||u_0||_{L^2_{(x+1)^mdx}}\\
&\leq& \frac{C_K\,C_1}{\sqrt{T}}||u_0||_{L^2_{(x+1)^mdx}}
\end{eqnarray*}
for $0<T\le 1$. Therefore 
\begin{equation}
\label{W5}
||u(t)||_{H^1_{(x+1)^m dx}} \le C_K\, C_1 e^\nu 
\frac{e^{-\nu t}}{\sqrt{t}} ||u_0||_{L^2_{(x+1)^mdx}}
\qquad \forall t\in (0,1).
\end{equation}
\eqref{W100} follows from \eqref{W5} and \eqref{L12}, since
$\mu < \nu$.
\qed
Let us return to the proof of Corollary \ref{c4}. Fix a number $\mu\in (0,\nu )$,
where $\nu$ is as in \eqref{L12}, and let us introduce the space
$$
F=\{u\in C(\R ^+; H^1_{ (x+1)^m dx}); \quad 
||e^{\mu t} u(t)||_{L^\infty (\R ^+; H^1_{(x+1)^mdx })} <\infty \}
$$
endowed with its natural norm. Note that \eqref{1} may be recast in the 
following integral form 
\begin{equation}
\label{W6}
u(t)=W(t)u_0 +\int_0^t W(t-s) N(u(s))\, ds
\end{equation}
where $N(u)=-uu_x$. We first show that \eqref{W6} has a solution in $F$
provided that 
$u_0\in H^1_{(x+1)^mdx}$ with $||u_0||_{H^1_{(x+1)^m  dx}}$ small enough.
Let $u_0\in H^1_{(x+1)^m dx}$ and $u\in F$ with 
$||u_0||_{H^1_{(x+1)^m dx}}\le r_0$ and $||u||_F\le R$, $r_0$ and $R$ 
being chosen later. We introduce the map $\Gamma$ defined by
\begin{equation}\label{gama}
(\Gamma u)(t)=W(t)u_0+\int_0^t W(t-s)N(u(s))\, ds\qquad  \forall t\ge 0.
\end{equation}
We shall prove that $\Gamma$ has a fixed point in the closed ball 
$B_R(0)\subset F$ provided that $r_0>0$ is small enough.

For the forcing problem
\begin{equation}
\begin{cases}
u_t + u_x + u_{xxx} +au = f\\
u(0,t) = 0 \\
u(x,0) = u_0(x)\nonumber
\end{cases}
\end{equation}
we have the following estimate 
\begin{eqnarray*}
&&\sup_{0\leq t\leq T}||u(t)||^2_{L^2_{(x+1)^m dx}} + 
\int_0^T \!\!\!\int_0^\infty (x+1)^{m-1} u_x^2 dxdt\\
&&\qquad \leq C\,\left(||u_0||^2_{L^2_{(x+1)^mdx }} + 
||f||^2_{L^1(0,T;L^2_{(x+1)^m dx}} \right) .
\end{eqnarray*}
Let us take $f=N(u)=-uu_x$. Observe that for all $x>0$
\begin{eqnarray*}
(x+1) u^2(x) 
& = & \left\vert \int_0^\infty \frac{d}{dx}[(x+1)u^2(x)]dx \right\vert \\
&\le&  C\left( 
\int_0^\infty (x+1)^m |u|^2 dx + 
\int_0^\infty (x+1)^{m-1} |u_x|^2 dx \right) 
\end{eqnarray*}
whenever $m\ge 2$. It follows that for some constant $K>0$
\begin{eqnarray*}
||uu_x||^2_{L^2_{(x+1)^m dx}} 
&\le& ||(x+1) u^2||_{L^\infty (\R ^+ )} \int_0^\infty (x+1)^{m-1}|u_x|^2dx\\
&\le& K ||u||^4_{H^1_{(x+1)^m dx}}.
\end{eqnarray*}
Therefore, for any $T>0$, 
\begin{eqnarray*}
&&\sup_{0\le t\le T} ||(\Gamma u)(t)||^2_{L^2_{(x+1)^m dx}}
+ \int_0^T\!\!\!\int_0^\infty (x+1)^{m-1}|(\Gamma u)_x|^2dxdt \nonumber\\
&&\qquad \le C\left(  ||u_0||^2_{L^2_{(x+1)^m dx }}  + 
\big(\int_0^T||u(t)||^2_{H^1_{(x+1)^mdx }}dt  \big)^2 \right) <\infty. 
\end{eqnarray*}
Thus $\Gamma u\in C(\R ^+, L^2_{(x+1)^mdx})\cap 
L^2_{loc}(\R ^+; H^1_{(x+1)^m dx})$ with $(\Gamma u)(0)=u_0$.
We claim that $\Gamma u\in F$. Indeed, by \eqref{L12}, 
$$
||e^{\mu t} W(t) u_0||_{H^1_{(x+1)^m dx}} \le C_1
||u_0||_{H^1_{(x+1)^m dx}}
$$
and for all $t\ge 0$
\begin{eqnarray*}
||e^{\mu t} \int_0^t W(t-s) N(u(s)) ds||_{H^1_{(x+1)^m dx}}
&\le& C e^{\mu t} \int_0^t
\frac{e^{-\mu (t-s)}}{\sqrt{t-s}} ||N(u(s))||_{L^2_{(x+1)^m dx}} ds\\ 
&\le& C\int_0^t \frac{e^{\mu s}}{\sqrt{t-s}}
K (e^{-\mu s}||u||_F)^2 ds \\
&\le& CK||u||^2_F \int_0^t\frac{e^{-\mu (t-s)}}{\sqrt{s}}ds\\
&\le& CK (2+\mu ^{-1}) ||u||^2_F
\end{eqnarray*}
where we used Lemma \ref{l2}.
Pick $R>0$ such that $CK(2 + \mu ^{-1})R \leq\frac{1}{2}$, and $r_0$ such that $C_1r_0 = \frac{R}{2}$. Then,
for $||u_0||_{H^1_{(x+1)^m dx}} \leq r_0$ and $||u||_F \leq R$, we obtain that
$$||e^{\mu t}(\Gamma u)(t)||_{H^1_{(x+1)^mdx}} \leq C_1r_0 + 
CK(2 + \mu ^{-1})R^2 \leq R, \quad t \geq 0.$$
Hence $\Gamma$ maps the ball $B_R(0)\subset F$ into itself. Similar computations show that
$\Gamma$ contracts. By the contraction mapping theorem, $\Gamma$ 
has a unique fixed point $u$ in $B_R(0)$. Thus 
$||u(t)||_{H^1_{(x+1)^mdx}} \le Ce^{-\mu t}||u_0||_{H^1_{(x+1)^m dx}}$
provided that $||u_0||_{H^1_{(x+1)^m dx}} \le r_0$ with
$r_0$ small enough. Proceeding as in the proof 
of Lemma \ref{l2}, we have that 
$$
||u(t)||_{H^1_{(x+1)^mdx}}
\le C\frac{e^{-\mu t}}{\sqrt{t}} ||u_0||_{L^2_{(x+1)^m dx}}
\qquad \mbox{ for } 0<t<1,$$
provided that $||u_0||_{L^2_{(x+1)^m dx}} \le \rho _0$ with $\rho _0<1$ 
small enough. The proof is complete with a decay rate $\mu' < \mu $. 
\qed
\begin{corollary}
\label{c5}
Assume that $a(x)$ satisfies \eqref{2} and that $\partial _x^ka\in L^\infty (\R ^+)$
for all $k\ge 0$. Pick any
$u_0\in L^2_{(x+1)^m dx}$. Then for all $\varepsilon >0$, all 
$T>\varepsilon$, and all $k\in \{1 ,... , m\}$, there exists a constant 
$C=C(\varepsilon,T, k)>0$ such that
\begin{eqnarray}
\int_\varepsilon^\infty (x+1)^{m-k} |\partial_x^k u(x,t)|^2dx 
\le C||u_0||^2_{L^2_{(x+1)^mdx}}\qquad \forall t\in [\varepsilon ,T]. 
\end{eqnarray}
\end{corollary}
\noindent {\bf Proof.}
The proof is very similar to the one in \cite[Lemma 5.1]{KF} and so we 
only point out the small changes.
First, it should be noticed that the presence in the KdV equation
of the extra terms $u_x$ 
and $a(x)u$ does not cause any serious trouble. On the other hand, choosing 
a cut-off function in $x$ of the form
$\eta(x)=\psi _0 (x/\varepsilon)$ (instead of $\eta (x)=\psi _0(x-x_0+2)$
as in \cite{KF})  
where $\psi _0\in C^\infty (\R , [0,1])$ satisfies $\psi_0(x)=0$ for
$x\le 1/2$ and $\psi_0(x)=1$ for $x\ge 1$, allows to overcome the fact
that $u$ is a solution of \eqref{1} on the half-line only.
\qed

\subsection{Decay in $L^2_b$}

This section is devoted to the exponential decay in
$L^2_b$. Our result reads as follows:

\begin{theorem}
\label{dec-exp}
Assume that the function $a=a(x)$
satisfies \eqref{2} with $4 b^3 + b < a_0$. 
Then, for all $R>0$,
there exist $C > 0$ and $\nu > 0$, such that
$$||u(t)||_{L^2_b}  \leq C  e^{-\nu t}||u_0||_{L^2_b}  \qquad t\ge 0$$
for any solution $u$ given by Theorem \ref{global-exp}.
\end{theorem}

\noindent {\bf Proof.} We introduce the Lyapunov function
\begin{equation}\label{exp1}
V(u) = \displaystyle\frac{1}{2}\int_0^\infty u^2 e^{2bx} dx +
c_b\int_0^\infty u^2 dx,
\end{equation}
where $c_b$ is a positive constant that will be chosen later. Then,
adding (\ref{fp30}) and (\ref{fp-4}) hand by hand we obtain
\begin{equation}\label{exp3}
\begin{array}{l}
\vspace{1mm}V(u) - V(u_0) = \displaystyle(4b^3 +
b)\int_0^T\int_{x_0}^\infty u^2 e^{2bx} dx dt
+ \displaystyle(4b^3 + b)\int_0^T\int_0^{x_0} u^2 e^{2bx} dx dt\\
\vspace{1mm}\qquad\qquad\qquad\qquad
-\,\,3b\,\displaystyle\int_0^\infty\int_0^\infty
u^2_x e^{2bx} dx dt +\frac{2b}{3}\int_0^T \int_0^\infty u^3
e^{2bx}dx
dt\\
\qquad\qquad\qquad\qquad- \,(c_b + \displaystyle\frac{1}{2})\int_0^T
u_x^2(0,t)dt - \,\, \displaystyle\int_0^T\int_0^\infty
a(x)|u|^2 (e^{2bx} +2c_b) dxdt,
\end{array}
\end{equation}
where $x_0$ is the number introduced in \eqref{2}. 
On the other hand, since $L^2_b \subset L^2_{(x+1)dx}$,
$||u(t)||_{L^2(0,\infty)}$ and $||u_x(t)||_{L^2(0,\infty)}$ decays
to zero exponentially. Consequently, from Moser estimate
we deduce that $||u(t)||_{L^\infty (0,\infty)}\rightarrow0$.
We may assume that $(2b/3)||u(t)||_{L^\infty} <\varepsilon =a_0-(4b^3+b)$
for all $t\ge 0$, by changing $u_0$ into $u(t_0)$ for $t_0$ large enough.
Therefore
\begin{equation}\label{exp3.1}
\begin{array}{l}
\vspace{1mm}\displaystyle\frac{2b}{3}\int_0^T\int_0^\infty |u|^3 e^{2bx}dx dt\\
\leq
\displaystyle\frac{2b}{3}\int_0^T ||u(t)||_{L^\infty(0,\infty)}
\left( \int_0^\infty |u|^2 e^{2bx}dx\right)  dt
\leq \varepsilon \int_0^T\int_0^\infty u^2 e^{2bx}dx dt.
\end{array}
\end{equation}
So, returning to (\ref{exp3}), the following
holds
\begin{equation}\label{exp4}
\begin{array}{l}
\vspace{1mm}V(u) - V(u_0)
-(4b^3+b +\varepsilon )\displaystyle\int_0^T\int_0^{x_0} u^2 e^{2bx} dx dt\\
+3b\displaystyle\int_0^T\int_0^\infty u^2_x e^{2bx} dx dt +
(c_b + \displaystyle\frac{1}{2})\int_0^T u_x^2(0,t)dt  +2 c_b
\displaystyle\int_0^T\int_0^\infty a(x)|u|^2dxdt \leq 0.
\end{array}
\end{equation}
Moreover, according to \cite{linares-pazoto} there exists $C > 0$ satisfying
\begin{equation}
\begin{array}{l}
\vspace{1mm}\displaystyle\int_0^T\int_0^{x_0} u^2 e^{2bx} dx dt\\
\leq
e^{2bx_0}\displaystyle\int_0^T\int_0^{x_0} u^2  dx dt  \leq C\,\{
\displaystyle\int_0^T u^2_x(0,t)dt +
\displaystyle\int_0^T\int_0^{\infty} a(x) u^2 dx dt \}\nonumber
\end{array}
\end{equation}
since $L^2_b \subset L^2(\R ^+)$. Then, choosing $c_b$ sufficiently large,
the above estimate and (\ref{exp4}) give us that
\begin{equation}\label{exp5}
\begin{array}{l}
\vspace{1mm}V(u) - V(u_0) \leq - C\,\{\displaystyle\int_0^T
u_x^2(0,t)dt + \displaystyle\int_0^T\int_0^\infty a(x) u^2
dx dt\\ \qquad\qquad\qquad\qquad + \displaystyle\int_0^T\int_0^\infty u^2_x e^{2bx} dx dt\}\leq - C\,V(u_0),
\end{array}
\end{equation}
which allows to conclude that $V(u)$ decays exponentially. The last
inequality is a consequence of the following results:

{\sc Claim 7.} There exists a positive constant $C>0$, such that
$$\int_0^T V(u(t)) dt \leq C \int_0^T\int_0^\infty u^2_x e^{2bx} dx dt.
$$

First, observe that
$$|\int_0^\infty u^2 e^{2bx} dx| = |-\frac{1}{b}\int_0^\infty uu_x e^{2bx}dx|
\leq \frac{1}{b} (\int_0^\infty u^2 e^{2bx}
dx)^{\frac{1}{2}}(\int_0^\infty u^2_x e^{2bx} dx)^{\frac{1}{2}},$$
therefore,
\begin{equation}
\begin{array}{l}\label{exp5.1}
\displaystyle\int_0^\infty u^2 e^{2bx} dx \leq
\frac{1}{b^2}\int_0^\infty u^2_x e^{2bx} dx.
\end{array}
\end{equation}
Then, from (\ref{2}) and (\ref{exp5.1}) we have
\begin{equation*}
V(u(t)) \le  (\frac{1}{2} + c_b) \int_0^\infty u^2 e^{2bx}dx 
\le (\frac{1}{2} + c_b)b^{-2} \int_0^\infty u_x^2 e^{2bx}dx 
\end{equation*}
which gives us Claim 7.

{\sc Claim 8.}
$$V(u_0) \leq C\,\{ \displaystyle\int_0^T u_x^2(0,t)dt +
\displaystyle\int_0^T\int_0^\infty u_x^2 e^{2bx} dxdt +
\displaystyle\int_0^T V(u(t)) dt\},$$ where $C$ is a positive constant.

Multiplying the first equation in (\ref{1}) by $(T-t)ue^{2bx}$ and
integrating by parts in $(0,\infty)\times (0,T)$, we obtain
\begin{equation}\label{exp7}
\begin{array}{l}
\vspace{1mm}-\displaystyle\frac{T}{2}\int_0^\infty |u_0(x)|^2
e^{2bx} dx + \frac{1}{2}\int_0^T\int_0^\infty |u|^2 e^{2bx} dx dt
+\,3b\displaystyle\int_0^T\int_0^\infty (T-t) u^2_x e^{2bx} dx dt\\
\vspace{1mm}+ \displaystyle\frac{1}{2} \int_0^T (T - t) u_x^2(0,t)dt
- \displaystyle(4b^3 + b)\int_0^T\int_0^\infty (T- t) u^2 e^{2bx} dx
dt\\
+\displaystyle\int_0^T\int_0^\infty (T - t) a(x)|u|^2 e^{2bx} dxdt
-\frac{2b}{3}\int_0^T \int_0^\infty (T - t) u^3 e^{2bx}dx dt= 0
\end{array}
\end{equation}
and therefore,
\begin{equation}\label{exp8}
\begin{array}{l}
\displaystyle \int_0^\infty |u_0(x)|^2 e^{2bx}
dx \leq C (\displaystyle\int_0^T u_x^2(0,t)dt +
\frac{1}{2}\int_0^T\int_0^\infty u^2 e^{2bx} dxdt
\\+\,\displaystyle\int_0^T\int_0^\infty u^2_x e^{2bx} dx dt
+\displaystyle\int_0^T
\int_0^\infty |u|^3 e^{2bx}dx dt).
\end{array}
\end{equation}
Then, combining (\ref{exp5.1}) and 
(\ref{exp3.1}), we derive Claim 8. \eqref{exp5} follows at once.
This proves the exponential decay when $||u(t)||_{L^\infty}\le 
3\varepsilon/(2b)$. The general case is obtained as in Theorem \ref{dec-pol}
\qed

\begin{corollary}
\label{c6} 
Assume that the function $a=a(x)$
satisfies \eqref{2} with $4 b^3 + b < a_0$. 
Then for any $R>0$, there exist positive constants 
$c=c(R)$ and $\mu = \mu (R)$ such that 
%and, in addition, $\frac{d^2 a}{d x^2}\in L^\infty$
\begin{equation}
||u_x(t)||_{L^2_b} \leq c 
\frac{e^{-\mu t}}{\sqrt{t}}
||u_0||_{L^2_b}
\end{equation}
for all $t>0$ and all $u_0\in L^2_b$ satisfying 
$||u_0||_{L^2_b}\le R$. 
\end{corollary}

\begin{corollary}\label{c7} 
Assume that the function $a=a(x)$
satisfies \eqref{2} with $4 b^3 + b < a_0$, and let $s\ge 2$. 
Then there exist some constants $\rho >0$, $C>0$ and $\mu >0$ such that  
%and, in addition, $\frac{d^2 a}{d x^2}\in L^\infty$
$$||u(t)||_{ H^s_b  } \leq C \frac{e^{-\mu t}}{t^{\frac{s}{2}}}
||u_0||_{ L^2_b }$$
for all $t>0$ and all $u_0\in L^2_b$ satisfying 
$|| u_0 ||_{L^2_b}\le \rho$. 
\end{corollary}
The proof of Corollary \ref{c6} (resp. \ref{c7}) is very similar to the proof
of Corollary \ref{c3} (resp. \ref{c4}), so it is omitted.

\section*{Acknowledgments.} This work was achieved while the first author
(AP) was visiting Universit\'e Paris-Sud with the support of the
Cooperation Agreement Brazil-France and the second author (LR) was visiting 
IMPA and UFRJ. 
LR was partially supported by the ``Agence Nationale de la Recherche'' (ANR), Project CISIFS, Grant ANR-09-BLAN-0213-02.

% You may incorporate your references as follows in your main tex file.
% Using BibTex is not recommended but can be handled.

%\medskip
% The data information below will be filled by AIMS editorial staff
%Received xxxx 20xx; revised xxxx 20xx.
%\medskip

\end{document}